\documentclass[11pt,letterpaper]{article}
\usepackage{amsmath}
\usepackage{amsfonts}
\usepackage{amssymb}
\usepackage{graphics}
\usepackage[dvips]{epsfig}
\newcommand{\RR}{\Bbb R}
\newcommand{\CC}{\Bbb C}
\newcommand{\PP}{\Bbb P}

\newcommand{\para}{\kern0.2em{\backslash} \kern-0.7em {\backslash} \kern0.2em }
\newcommand{\parap}{          {\backslash} \kern-0.7em {\backslash} \kern0.2em }

\newcommand{\thm}{Theorem }
\newcommand{\mfd}{manifold }

\newcommand{\subm}{submanifold }
\newcommand{\subms}{submanifolds }

\newcommand{\Kasten}{\hfill $\sqcap\!\!\!\!\sqcup$\\}
\newcommand{\nbhd}{neighborhood }

\newcommand{\riem}{Riemannian }
\newcommand{\tn}{tubular neighborhood }

\newcommand{\vs}{vector space }
\newcommand{\vb}{vector bundle }

\newcommand{\rem}{\textbf{Remark:} }
\newcommand{\e}{\epsilon}
\newcommand{\iso}{isotropic }
\newcommand{\Ws}{Weinstein's }

\newcommand{\av}{average }
\newcommand{\tnng}{\exp_{N_g} (\nu N_g)_1}

\newcommand{\tnngl}{\exp_{N_g} (\nu N_g)_{L}}
\newcommand{\tnnl}{\exp_N (\vertg )_{L}}

\newcommand{\vg}{{\varphi_g}}
\newcommand{\vgi}{{\varphi_g^{-1}}}
\newcommand{\og}{\omega_g}
\newcommand{\iog}{\int_g \omega_g}
\newcommand{\ot}{\omega_t}
\newcommand{\vfs}{vector fields }
\newcommand{\vf}{vector field }
\newcommand{\wrt}{with respect to }
\newcommand{\png}{\pi_{N_g}}
\newcommand{\vertg}{{Vert}^g}
\newcommand{\horg}{{Hor}^g}
\newcommand{\avg}{{aVert}^g}
\newcommand{\ahg}{{aHor}^g}
\newcommand{\lc}{{LC}^g}
\newcommand{\pf}{\textit{Proof: }}
\newcommand{\jf}{Jacobi-field }
\newcommand{\jfs}{Jacobi-fields }
\newcommand{\st}{such that }
\newcommand{\tj}{\tilde{J}}
\newcommand{\la}{\langle}
\newcommand{\ra}{\rangle}

\newcommand{\tr}{triangle }
\newcommand{\trs}{triangles }
\newcommand{\geo}{geodesic }

\newcommand{\pt}{point }

\newcommand{\tb}{\tilde{B}}
\newcommand{\tc}{\tilde{C}}
\newcommand{\tbc}{Q}
\newcommand{\tsigma}{\tilde{\sigma}}
\newcommand{\dsigma}{\dot{\sigma}}
\newcommand{\dtsigma}{\dot{\tilde{\sigma}}}
\newcommand{\dgamma}{\dot{\gamma}}
\newcommand{\tgamma}{\tilde{\gamma}}

\newcommand{\gammai}{{\gamma}^{b}}
\newcommand{\expi}{{\exp}^{-1}}
\newcommand{\pii}{{\pi}^{b}}
\newcommand{\np}{{\nabla}^{\perp}}
\newcommand{\ndt}{\frac{\nabla}{dt}}
\newcommand{\npdt}{\frac{\np}{dt}}

\newcommand{\npdto}{\frac{\np}{dt}}
\newcommand{\prt}{{pr}_{\dgamma(t)}}
\newcommand{\pro}{{pr}_{\dgamma(0)}}
\newcommand{\re}{R^{\epsilon}}

\newcommand{\te}{\text{\textit{tub}}^{\epsilon}}
\newcommand{\td}{\tilde{\delta}}
\newcommand{\de}{D^{\epsilon}}

\newcommand{\ag}{\alpha^g}
\newcommand{\bg}{\beta^g}
\newcommand{\qng}{Q_{N_g}}
\newcommand{\qgn}{Q_N^g}
\newcommand{\pgn}{\pi_N^g}
\newcommand{\rngt}{{(\rho_{N_g})_t}}
\newcommand{\rgnt}{{(\rho_{N}^g)_t}}
\newcommand{\onb}{orthonormal basis }

\newcommand{\orth}{orthogonal }
\newcommand{\otr}{orthogonal-to-radial }
\newcommand{\pr}{{pr}}
\newcommand{\npds}{\frac{\np}{ds}}
\newcommand{\nds}{\frac{\nabla}{ds}}
\newcommand{\lep}{L^{\e}}
\newcommand{\g}{\gamma}
\newcommand{\gd}{\dot{\gamma}}

\newcommand{\cg}{\mathfrak{g}}
\begin{document}
\ifx\href\undefined\else\hypersetup{linktocpage=true}\fi

\newtheorem{Thm}{Theorem}
\newtheorem{Lem}{Lemma}[section]
\newtheorem{Prop}{Proposition}[section]
\newtheorem{Cor}{Corollary}[section]

\newcommand {\comment}[1]{{\marginpar{*}\scriptsize{\bf Comments:}\scriptsize{\ #1 \ }}}

\title {Submanifold averaging in riemannian and symplectic
geometry}
\author{\begin{Large}Marco Zambon\end{Large}\\
Department of Mathematics, UC Berkeley, CA, USA\\
\textsc{zambon@math.berkeley.edu}}
\date{\today}
\maketitle

\begin{abstract} We give a construction to obtain
canonically an ``isotropic average'' of given $C^1$-close
 isotropic submanifolds of a symplectic manifold.
To do so we use an improvement of Weinstein's submanifold
averaging theorem (obtained in collaboration with H. Karcher) and
apply ``Moser's trick''. We also present an application to
Hamiltonian group actions.
\end{abstract}

\tableofcontents

\section {Introduction}\label{s1}

\noindent In 1999 Alan Weinstein [We] presented a procedure to
average a family $\{N_g\}$ of submanifolds of a Riemannian
manifold $M$: if the submanifolds are close to each other in a
$C^{1}$ sense, one can produce \emph{canonically} \footnote{The
construction is canonical because it does not involve any
arbitrary choice but uses only the Riemannian metric on $M$.}
an ``average'' $N$ which is close to each member of the family $\{N_g\}$. The main property of this averaging procedure is that it is equivariant with respect to isometries of $M$, and therefore if the family  $\{N_g\}$ is obtained by applying the isometric action of a compact group $G$ to some submanifold $N_0$ of $M$, the resulting average will be invariant under the $G$-action. This generalizes results about fixed points of group actions [We].\\
 In the first part of this paper we will exhibit a result by Hermann Karcher and the author
which improves Weinstein's theorem.\\[0.4cm]
In the main body of the paper we specialize Weinstein's averaging
to the setting of symplectic geometry: given a family of
\emph{isotropic} submanifolds  $\{N_g\}$ of a symplectic  manifold
$M$, we obtain an \emph{isotropic}
 average $L$.
We achieve this in two steps: first we introduce a compatible
Riemannian metric on $M$ and apply Weinstein's averaging to obtain
a submanifold $N$. This submanifold will be ``nearly isotropic''
because it is $C^1$-close to isotropic ones, and using the family
$\{N_g\}$ we will deform $N$ to an isotropic submanifold
$L$.\footnote{It would be interesting to find a way to deform any
given ``nearly isotropic'' submanifold  to an honest isotropic one
in a canonical fashion.} Our construction depends only on the the
symplectic structure of $M$ and on the choice of compatible
metric. Therefore applying our construction to the case of compact
group actions by isometric
 symplectomorphism we can obtain isotropic submanifolds which are invariant under the action.\\[0.4cm]
As a simple application we show that the image of an almost
invariant isotropic submanifold
under a compact Hamiltonian action is ``small''.\\[0.4cm]
Another application is the following: given a symplectic action of
a compact group $G$ on two symplectic manifolds $M_1$ and $M_2$
together with an almost equivariant symplectomorphism $\phi:M_1
\rightarrow M_2$, apply the averaging procedure to
$graph(\phi)\subset M_1\times M_2$. If the resulting $G$-invariant
submanifold $L$ is a graph, then it will be the graph of a
$G$-equivariant symplectomorphism. This means that we would be
able to deform almost equivariant symplectomorphisms to
equivariant ones. To assure that $L$ is again a graph one needs to
improve Weinstein's averaging procedure \footnote{We need to
improve Weinstein's theorem in order to assure that $graph(\phi)$
and $L$ be $C^1$-close; see Remark 1 in Section \ref{s8}}
; this is the subject of work in progress.\\[0.4cm]
We would like to extend our averaging procedure to coisotropic submanifolds too: indeed,
 if one could average any two coisotropic submanifolds $N_0$ and $N_1$
which are close to each other, then by ``shifting weights'' in the parameter space $G=\{0,1\}$ one would
produce a continuous path of coisotropic submanifolds connecting $N_0$ to $N_1$.
This would show that the
space of coisotropic submanifolds is locally path connected.\\[0.4cm]
In the remainder of the introduction we will recall the averaging
procedure in the Riemannian setting by Weinstein (see [We]), we
will state our results, and we will outline our construction to
average isotropic submanifolds.

\subsection{Averaging of \riem \subms}\label{s11}

The starting point for our isotropic averaging construction is the
statement of Theorem 2.3 in [We]. We first recall some definitions
from [We] in
order to state the theorem.\\[0.4cm]
If $M$ is a Riemannian \mfd and $N$ a submanifold,
$(M,N)$ is called \textit{gentle pair} if (i) the
normal injectivity radius of $N$ is at least one; (ii)
the sectional curvatures of $M$ in the tubular \nbhd
of radius one about $N$ are bounded in absolute value
by one; (iii) the injectivity radius of each point of
the above \nbhd is at least one.\\[0.4cm]
The \textit{distance between two subspaces} of the
same dimension $F,F'$ of a Euclidean vector space $E$,
denoted by $d(F,F')$, is equal to the $C^0$-distance
 between the unit spheres of
$F$ and $F'$ considered as \riem \subms of the unit
sphere of $E$. This distance is symmetric and
satisfies $d(F,F')=d(F^{\perp},{F'}^{\perp})$. It is
always smaller or equal than $\frac{\pi}{2}$, and it
is equal to $\frac{\pi}{2}$ iff $F$ and ${F'}^{\perp}$
are not transversal.\\[0.4cm]
One can define a \textit{$C^1$-distance between two
\subms} $N, N'$ of a \riem \mfd if $N'$ lies in the
\tn of $N$ and is the image under the normal
exponential map of $N$ of a section of $\nu
N$ (so $N$ and $N'$ are necessarily diffeomorphic).
This is done by assigning two numbers to each $x'
\in N'$: the length of the geodesic segment from $x'$
to the nearest point $x$ in $N$ and the distance
between $T_{x'}N'$ and the parallel translate of
$T_xN$ along the above geodesic segment. The
$C^1$-distance is defined as the supremum of these
numbers as $x'$ ranges over $N'$ and is denoted by
$d_1(N, N')$.\\
Note that this distance is not symmetric, but if $(M,N)$ and $(M,
N')$ are both gentle pairs with $d_1(N, N')<\frac{1}{4}$, then
$d_1(N', N)<250 d_1(N, N')$ (see Remark 3.18 in [We]).\\[0.4cm]
The improvement of Theorem 2.3 in [We] by Karcher and the author is our Theorem 4 and reads \footnote{We omit the
compactness assumption on the $N_g$'s stated there since it is
superfluous.}:
\\[0.2cm]
\textbf{Theorem [Weinstein] } \textit{Let $M$ be a \riem \mfd and
$\{N_g\}$ a family of submanifolds of $M$ parametrized in a
measurable way by elements of a probability space $G$, such that
all the pairs $(M,N_g)$ are gentle. If
$d_1(N_g,N_h)<\epsilon<\frac{1}{20000}$ for all $g$ and $h$ in
$G$, there is a well defined \emph{\textbf{center of mass}}
submanifold $N$ with $d_1(N_g,N)<2500 \epsilon$ for all $g$ in $G$. The
center of mass construction is equivariant with respect to
isometries of $M$ and measure preserving
automorphisms of $G$.}\\[0.4cm]
\rem For any $g\in G$ the center of mass $N$ is the image under the exponential map of a section of $\nu N_g$ and $d_0(N_g,N)<100\e$.\\[0.4cm]
From this one gets immediately a statement about invariant \subms
under compact group actions (cfr. \thm
2.2 of [We]).\\[0.4cm]

\subsection{Averaging of isotropic \subms}\label{s12}

Recall that for any symplectic manifold $(M,\omega)$ we can choose
a compatible Riemannian metric $g$, i.e. a metric such that the
endomorphism $I$ of $TM$ determined by $\omega(\cdot\;,
I\cdot)=g(\cdot\;,\;\cdot)$ satisfies $I^2=-Id_{TM}$. The tuple
$(M,g,\omega, I)$ is called \emph{almost-K\"ahler manifold}. To
prove our Main Theorem we need assume a bound on the $C^0$-norm of
$\nabla \omega$ (here $\nabla$ is the Levi-Civita connection given
by $g$), which measures how far our almost-K\"ahler manifold is
from being K\"ahler \footnote{Recall that an almost-K\"ahler
manifold is K\"ahler if the almost complex structure $I$
is integrable, or equivalently if $\nabla I=0$ or $\nabla \omega=0$.}.
We state the theorem choosing the bound to be one (but see Remark i)).\\
\begin{Thm}\emph{[Main Theorem]} Let $(M^m,g,\omega,I)$ be an almost-K\"ahler \mfd
satisfying $|\nabla \omega|<1$
and $\{N^n_g\}$ a family of isotropic submanifolds of $M$
parametrized in a measurable way by elements of a probability
space $G$, such that all the pairs $(M,N_g)$ are gentle. If
$d_1(N_g,N_h)<\epsilon<\frac{1}{70000}$ for all $g$ and $h$
in $G$, there is a well defined \emph{\textbf{isotropic center of
mass}} submanifold $L^n$ with $d_0(N_g,L)<1000 \e$ for all $g$
in $G$. This construction is equivariant with respect to
isometric symplectomorphisms of $M$ and measure preserving
automorphisms of $G$.\end{Thm}
 \rem
i) The theorem still holds if we assume higher bounds on $|\nabla
\omega|$, but in this case the bound $\frac{1}{70000}$ for $\epsilon$
would have to be chosen smaller. See the remark in Section \ref{s74}.\\
ii) Notice that we are not longer able to give estimates on the
$C^1$-distance of the \iso center of mass from the $N_g$'s. Such
an estimate could possibly be given provided we have more
information about the extrinsic geometry of \Ws center of mass
submanifold; see Remark 1 in Section \ref{s8}. Instead we can only
give estimates on the  $C^0$-distances
$d_0(N_g,L)= {\sup}\{d(x,N_g):x\in L \}$.\\[0.4cm]

An easy consequence of our Main Theorem is a statement about group
actions. Recall that, given any action of a compact Lie group $G$
on a symplectic manifold $(M,\omega)$ by symplectomorphisms, by
averaging over the compact group one can always find some
invariant metric $\tilde{g}$. Using $\omega$ and  $\tilde{g}$ one
can canonically construct a metric $g$ which is compatible with
$\omega$ (see [Ca]), and since $g$ is constructed canonically out
of objects that are $G$-invariant, it will be $G$-invariant too.
Therefore the group $G$ acts respecting the structure of the
almost K\"ahler manifold $(M,g,\omega)$. In general it does not
seem possible to give any bound on $|\nabla \omega|$, where
$\nabla$ is the Levi-Civita connection corresponding to $g$.
\begin{Thm} Let  $(M,g,\omega,I)$ be an almost-K\"ahler
\mfd satisfying $|\nabla \omega|<1$ and let $G$ be a compact Lie
group acting on $M$ by isometric symplectomorphisms. Let ${N_0}$
be an \iso \subm of $M$ such that $(M,N_0)$ is a gentle pair and
$d_1({N_0},g{N_0})<\e<\frac{1}{70000}$ for all $g \in G$. Then
there is a $G$-invariant \iso \subm $L$ with $d_0({N_0},L)<
1000\e$. \end{Thm}

 The invariant isotropic submanifold $L$ as
above is constructed endowing $G$ with the bi-invariant
probability measure and applying Theorem 1 to the family
$\{gN_0\}_{g\in G}$. The resulting isotropic average $L$ is
$G$-invariant because of the equivariance properties of the
averaging procedure.
\subsection{Outline of the proof of the Main Theorem}\label{s13}

This is the main subsection of this paper.\\
 We will try to convince the reader that  the
construction we use to prove Theorem 1 works if only one chooses
$\e$ small
enough. Let us begin by requiring $\e <\frac{1}{20000}
$.\\[0.4cm]
$\bullet$ \textbf{PART I } We start by considering the \av of the
\subms $N_g$ as in \thm 2.3 of [We], which we will denote $N$. We
will use the notation $\exp_N$ to indicate the restriction of the
exponential map to $TM|_N$, and similarly for any of the $N_g$'s. For any
$g$ in $G$, the \av $N$ lies in a \tn of $N_g$ and is the image
under $\exp_{N_g}$ of a section $\sigma$ of $\nu N_g$ (see [We]).
Therefore for any point $p$ of $N_g$ there is a canonical path
$\gamma_q(t)=\exp_p(t\cdot \sigma(p))$ from $p$ to the unique
point $q$ of $N$ lying in the normal slice of $N_g$ through $p$.
Here, using the notation $ (\nu N_g)_1$ for the open unit disk
bundle in $\nu N_g$, we denote by the term ``normal slice'' the
submanifold $\exp_{N_g}(\nu_p N_g)_1$. We define the following
map:
$$\varphi_g: \tnng \rightarrow M, \;\;\; \exp_p(v)
\longmapsto \exp_q(_{\gamma_q}\parap v).$$
Here $p,q$, and $\gamma_q$ are as above,  $v\in
(\nu_p  N_g)_1$.\\
So $\varphi_g$ takes the normal slice $\exp_p( \nu_p
N_g)_1$ to $\exp_q(Vert_q^g)_1$, where $Vert_q^g
\subset T_qM$ is the parallel translation along
$\gamma_q$ of $\nu_pN_g \subset T_pM$.\\
We have $d(Vert_q^g,\nu_qN) <d_1(N_g,N)<2500 \e < \frac{\pi}{2}$,
so $Vert_q^g$ and $T_qN$ are transversal. Therefore $\vg$ is a
local diffeomorphism at all points of $N_g$, and it is clearly
injective there. Using the geometry if $N_g$, $N$ and $M$ in
Proposition \ref{Prop51} we will show that $\vg$ is a
diffeomorphism onto if restricted to the tubular \nbhd $\exp_{N_g}
(\nu N_g)_{0.05}$ of $N_g$.\\ We restrict our map to this
neighborhood and we also restrict the target space so to obtain a
diffeomorphism, which we will
still denote by $\vg$ .\\[0.4cm]

\centerline{\epsfig{file=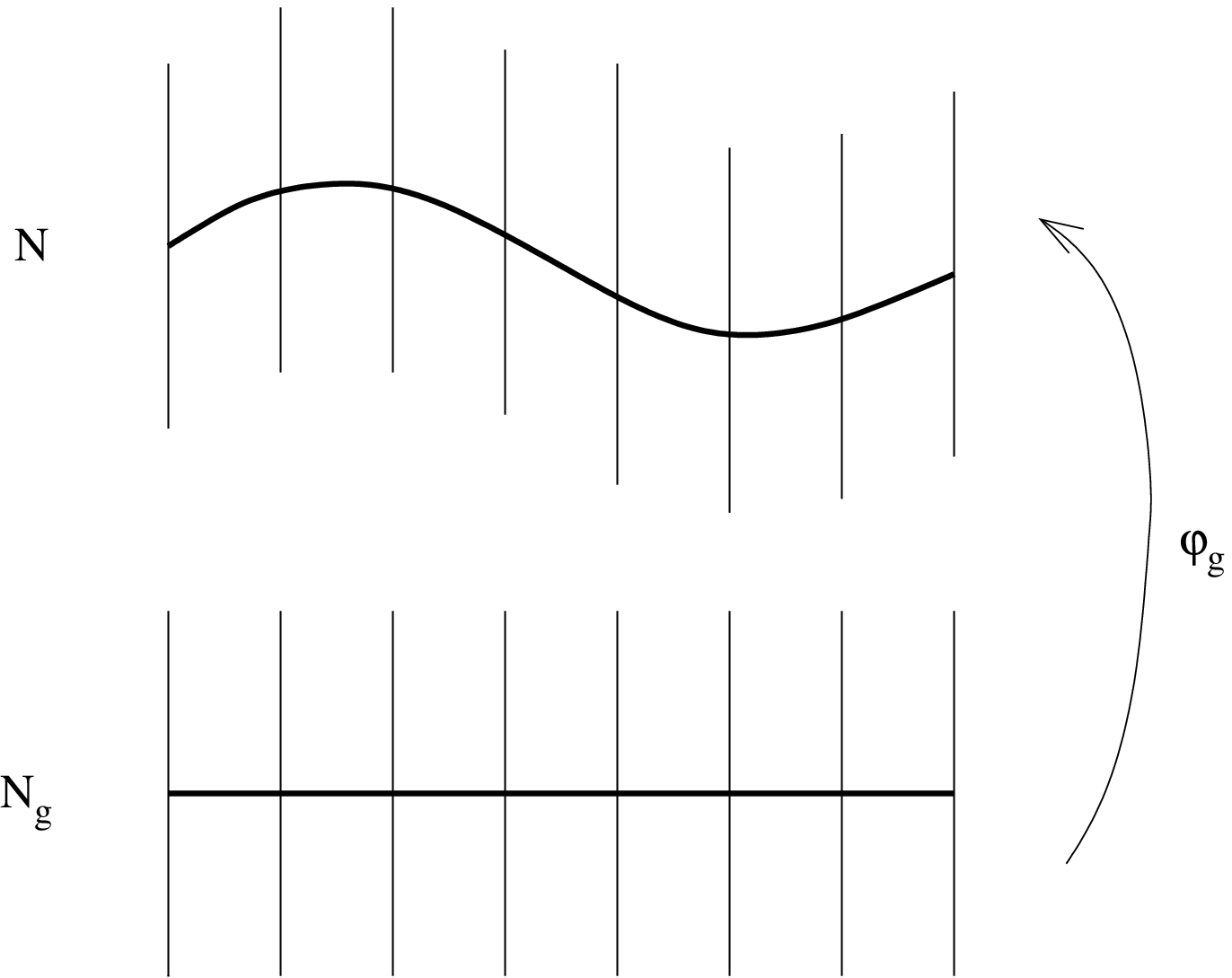,height=7cm}}

\noindent $\bullet$ \textbf{PART II }Now we introduce the
symplectic form
$$\omega_g:=(\vgi)^* \omega$$ on  $\exp_N (Vert^g)_{0.05}$. Notice
that $N$ is isotropic with respect to $\omega_g$ by construction,
therefore it is also with respect to the two-form $\iog$ which is
defined on $\cap_{g\in G}\exp_N(\vertg)_{0.05}$. We would like to
apply Moser's trick \footnote{ Recall that Moser's Theorem states
the following: if $\Omega_t$ ($t\in[0,1]$) is a smooth family of
symplectic forms lying in the same cohomology class in a compact
manifold then there is a family of diffeomorphisms $\rho_t$ with
$\rho_0=Id$ satisfying $\rho_t^*\Omega_t=\Omega_0$.} (see [Ca,
Chapter III]) to $\omega$ and $\iog$. To do so first we restrict
our forms to a smaller tubular  neighborhood $\te$ of $N$, which we define in Section \ref{s71}. To
apply Moser's trick
we have to check:\\
\begin{enumerate}\item  \textit{ On $\te$ the convex linear
combination $\ot=\omega +t(\iog -\omega)$ consists of symplectic
forms.}\\
Indeed we will show that on $\te$ the differential of ${\vgi}$ is
``close'' to the parallel translation $\,\para$ along certain
``canonical'' geodesics that will be specified at the beginning of
Section \ref{s3}. This and the bound on $|\nabla \omega|$ imply
that for any $q\in \te$ and nonzero $X,Y\in T_qM$:
\begin{equation*}\begin{split}(\og)_q(X,Y)=&
\,\omega_{\vgi(q)}\left(\vgi_*(X), \vgi_*(Y) \right)
\\ \approx &
\,\omega_{\vgi(q)}\Big(\para X,
\para Y \Big)\\ \approx & \,\omega_q(X,Y),
\end{split}\end{equation*}
i.e. $\omega_g$ and $\omega$ are very close to each other.
 So $\omega_t(X,IX)\approx \omega(X,IX)=|X|^2>0$. Therefore each $\omega_t$
is non-degenerate, and it is clearly also closed. \\
 \item \textit{On $\te$ the forms $\omega$ and $\iog$ belong to the
same cohomology class}.\\
Fix $g\in G$. The inclusion $i: \te \hookrightarrow \exp_{N_g}(\nu
N_g)_{1}$ is homotopic to $\vgi:\te
  \rightarrow
   \exp_{N_g}(\nu N_g)_{1}$. A homotopy is given by thinking
of $N$ as a section of $\nu N_g$ and ``sliding along the fibers''
to the zero-section. Therefore these two maps induce the same map
in cohomology, and pulling back $\omega$ we have
$$[\omega|_{\te}]=i^*[\omega]=(\vgi)^*[\omega]=[\omega_g].$$  Integrating over $G$ finishes
the argument.\\[0.2cm]
\end{enumerate}
Now we can apply Moser's trick: if $\alpha$ is a one
form on $\te$ such that $d \alpha$ is equal to
$\frac{d}{dt} \ot = \iog -\omega$, then the flow
$\rho_t$ of the time-dependent \vf
$v_t:=-\tilde{\ot}^{-1}(\alpha)$ has the property
$\rho_t^* \ot =\omega$ (and in particular
$\rho_1^*(\iog)=\omega$) where it is defined.
Therefore if $L:=\rho_1^{-1}(N)$ is a well-defined
submanifold of $\te$, then it will be \iso \wrt
$\omega$ since $N$ is \iso \wrt $\iog$.\\
We will construct canonically a primitive $\alpha$ as above in
Section \ref{s72}. Using the fact that the distance between the
$N_g$'s and $N$ is small, we will show that $\alpha$ has small
maximum-norm. So, if $\e$ is small enough, the time-1 flow of the
time-dependent \vf $\{ -v_{1-t} \}$ will not take $N$ out of $\te$
and
$L$ will be well-defined.\\[0.4cm]
Since our construction is canonical after fixing the
almost-K\"ahler structure $(g, \omega, I)$ of $M$ and the
probability space $G$, the construction of $L$ is
equivariant \wrt isometric symplectomorphisms of $M$ and
measure preserving automorphisms of $G$.\\

\subsection{Structure of the paper and acknowledgements}\label{s14}
This paper is organized as follows:\\
In Section 2 we present the improvement of Theorem 2.3 in
[We] obtained by Karcher and the author. In Section 3 we will start the proof of the Main Theorem by
studying the map $\vg$. In Section 4 we will state a proposition
about geodesic triangles, and in Section 5 we will apply it in our
setup. This will allow us to show in Section 6 that each $\vg$ is
injective on $\exp_{N_g}(\nu N_g)_{0.05}$. The proofs of some
estimates of Sections 4 and 5 are rather involved, and we present
them in the three appendices. This will conclude the
 proof of the first part of the theorem.\\
In Section 7 we will make use for the first time of the symplectic
structure of $M$. We will show that the $\ot$'s are symplectic
forms and that the 1-form $\alpha$, and therefore the Moser \vf
$v_t$, are small in the maximum norm. Comparison with the results
of Section 6 will end the proof of the Main
Theorem.\\
Section 8 will be devoted to
remarks about the Main Theorem, and in Section 9 we will present a simple
 application to Hamiltonian group actions.\\[0.4cm]
At this point I would like to thank everyone who helped me and
supported me in the preparation of this paper. In particular I
would like to thank Alan Weinstein for helpful discussions during
the preparation of this paper, the referee for his careful review
of the manuscript, his interest and for suggesting improvements,
Yael Karshon for proposing the application in Section 9 and River
Chiang for simplifying the arguments used there. Further I thank
Hermann Karcher for sharing the ideas involved in Section 2 and for the collaboration.

\section{Improved error estimates for the shape operators of parallel tubes \-
with application to Weinstein's submanifold averaging}\label{s2}

In this section we will present the improvement of
Theorem 2.3 in [We] obtained by Hermann Karcher and the author. In the first subsection we will improve
Proposition 3.11 in [We]. Then using this result we will follow
Weinstein's proof and present the statement of the
improved theorem.

\subsection{Estimates for the shape operators of parallel tubes }\label{s21}

In Proposition 3.11 of [We] one has the setup we are going to
describe now.\\
$M$ is a Riemannian manifold, $N$ is a submanifold of $M$ such
that $(M,N)$ form a gentle pair (so the second fundamental form
$B$ of $N$ satisfies $|B|\le \frac{3}{2}$, see [We, Cor 3.2]). In
the tubular neighborhood of radius one about $N$ let $\rho_N$ be
the distance function from $N$, and $P_N=\frac{1}{2}\rho^2_N$. We
are interested in estimating the Hessian of $P_N$, i.e. the
symmetric endomorphism of each tangent space of the tubular
neighborhood given by $H_N(v)= \nabla_v \text{grad} P_N.$
Differentiating the relation $ \text{grad} P_N=\rho_N\cdot
\text{grad} \rho_N$ we see that
$$H_N(v)= \langle U_N,v \rangle U_N + \rho_N\cdot S_N(\text{pr}(v))$$
where $U_N=\text{grad} \rho_N$ is the radial unit vector (pointing
away from $N$), $\text{pr}$ denotes orthogonal projection onto
$U_N^{\perp}$, and $S_N$ is the second fundamental form
 of the tube given by a level
set $\tau(t)$ of $\rho_N$ in direction of the normal vector
$U_N$\footnote{ So $S_Nv=\text{pr}(\nabla_v U_N)$ for all vectors
$v$ tangent to $\tau(t)$, where $\nabla$ is the Levi-Civita
connection on $M$.}
.\\
Proposition 3.11 of [We] states that, at a point $p$ of distance
$t\le \frac{1}{4}$ from $N$, the following estimate holds for the
decompositions into vertical and horizontal parts \footnote{See
our Section \ref{s3} or Section 2.1 in [We] for the definition of
vertical and horizontal bundle at $p$.}
 of $T_pM$:
$$\begin{bmatrix}   0.64\cdot I&0 \\0 & -3t \cdot I
\end{bmatrix} <
H_N <
\begin{bmatrix}  1.32 \cdot I&0 \\0 &3t\cdot I
\end{bmatrix},$$
where for two symmetric matrices $P$ and $Q$ the inequality $P<Q$
means that $Q-P$ is positive definite.\\
The above proposition is proved using the Riccati equation. An
immediate consequence is Corollary 3.13 in [We], which states that, if $v$
is a horizontal vector and $w$ a vertical vector at $p$, then $|\langle H_N(v), w
\rangle| \le 3 \sqrt{t}|v||w|.$ This square root is
responsible for the presence of upper bounds proportional to
$\sqrt{\epsilon}$ rather than $\epsilon$ in Theorems 2.2 and 2.3 of
[We].\\
\noindent We will improve the estimate of Corollary 3.13 in [We]
determining $S_N$ by means of Jacobi-fields estimates rather than
by the Riccati equation. More precisely, we will make use of
this simple observation:

\begin{Lem}\label{Lemk1}
Let $N$ be a submanifold of the Riemannian manifold $M$, and fix
$t\le \text{normal injectivity radius of }N$. Let p lie in the
tube $\tau(t):=\rho_N^{-1}(t)$, and let $S_N:T_p\tau(t)
\rightarrow T_p\tau(t)$ be the second fundamental form in
direction of $U_N$. For any $v\in T_p \tau(t)$ consider the Jacobi
field $\tilde{J}(r)$ arising from the variation $r \mapsto
\exp_{c(s)}rU_N(c(s))$, where $c(s)$ is any curve in $\tau(t)$
tangent to $v$. Then
$$S_Nv=\tilde{J}'(0).$$
\end{Lem}
\textit{Proof: } Denoting by $f(s,r)$ the above variation and by
$\nabla$ the Levi-Civita connection on $M$ we have
\[ \begin{split}
\tilde{J}'(0)&=\frac{\nabla}{dr}|_0 \frac{d}{ds}|_0 f(s,r)\\
&=\frac{\nabla}{ds}|_0 \frac{d}{dr}|_0 f(s,r)\\
&=\frac{\nabla}{ds}|_0      U_N(c(s))\\
&=\nabla_v U_N\\
&=\text{pr} (\nabla_v U_N)\\
&=S_N v.
\end{split}
\]
\Kasten

 \noindent
Using the above lemma we will be able to prove this improvement of
Proposition 3.11 in [We], for which we don't require $(M,N)$ to be
a gentle pair but only a bound on $|B|$:

\begin{Thm}
Let $N$ be a submanifold of the Riemannian manifold $M$ with
second fundamental form $B$, and fix $t\le \text{normal
injectivity radius of }N$. Let $\gamma$ be a unit speed geodesic
emanating normally from $N$. Let $\tau(t)$ be the $t$-tube about
$N$, and let $S_N(t)$ denote the second fundamental form of
$\tau(t)$ in direction $\dot{\gamma}(t)$ at $\g(t)$. Then w.r.t.
the splitting into vertical and horizontal spaces of
$T_{\gamma(t)}\tau(t)$, as long as $t\le
\min\{\frac{1}{2},\frac{1}{2|B|} \}$, we have
$$t\cdot S_N(t)=
\begin{bmatrix}   I&0 \\0 &tB
\end{bmatrix}+
\begin{pmatrix} 16 t^2 & 16 t^2\\ 16 t^2 & (22+2|B|^2 )t^2
\end{pmatrix}.
$$
\end{Thm}
\textbf{Remark :} We adopt the following unconventional notation:
If $M,\tilde{M}$ are matrices and $c$ a real number,
$M=\tilde{M}+c$ means that $M-\tilde{M}$ has operator norm $\le
c$. Generalizing to the case where we consider also
vertical-horizontal decompositions of matrices,
$$\begin{bmatrix}  A&B \\C & D
\end{bmatrix}\le \begin{bmatrix}  \tilde{A}&\tilde{B} \\\tilde{C} &
\tilde{D}
\end{bmatrix}+\begin{pmatrix}  a&b \\c & d
\end{pmatrix}$$ means that the above convention holds for each
endomorphism between horizontal/vertical spaces, i.e.
$A-\tilde{A}$ has operator norm $\le a$ and so on.\\

\noindent
 \textit{Proof: }  Choose an orthonormal
basis $\{E_1,\dots,E_{n-1}\}$ of $\gd(0)^{\perp}\subset
T_{\gamma(0)}M$ such that $\{E_1,\dots,E_k\}$ lie in the normal
space to $N$ and $\{E_{k+1},\dots,E_{n-1}\}$ lie in the tangent
space to $N$. (Here $\dim(M)=n$.) Now we define Jacobi fields
$J_i$ along $\g$ with the following initial conditions:
\begin{equation*}
\begin{cases}
J_i(0)=0,\,\,J_i'(0)=E_i & \text{ if $i\le k$
(vertical Jacobi fields)}\\
 J_i(0)=E_i,\,\,J_i'(0)=B_{\gd(0)}E_i  & \text{ if $i\ge k+1$  (horizontal
 Jacobi fields).}\\
\end{cases}
\end{equation*}
Notice that, among all $N$-Jacobi fields (see Section \ref{s3} for
their definition) satisfying $J_i(0)=E_i$, our $J_i$ are those
having smallest derivative at time zero. Also notice that all
$J_i$ and their derivatives are perpendicular to $\gd(0)$,
therefore, as long as the $J_i(t)$ are linearly independent, they
form a basis of $\gd(t)^{\perp}=T_{\g(t)}\tau(t)$. Also, the
$J_i$'s are $N$-Jacobi fields, i.e. Jacobi fields for which
$J_i(0)$ is tangent to $N$ and $J'_i(0)-B_{\gd(0)}J_i(0)$ is
normal to $N$, or equivalently Jacobi fields that arise from
variations of geodesics emanating normally from $N$ (see [Wa, p.
342]). Moreover the $J_i$'s are a basis of the space of $N$-Jacobi
fields along $\gamma$ which are orthogonal to $\gd$, and this
space coincides with the space of $N$-Jacobi fields  arising from
a variation of unit-speed \footnote{Indeed, connecting the points
of an integral curve in $\tau(t)$ of some $J_i(t)$ to $N$ by unit
speed shortest geodesics we obtain such a variation, and the
Jacobi field arising from this variation must be $J_i$ since it is
an $N$-Jacobi field orthogonal to $\gd$ which coincides with
$J_i(t)$ at time $t$.} geodesics normal to $N$. The velocity
vectors of such variations at time $t$ coincide with $U_N$.
Therefore applying Lemma \ref{Lemk1} with $v=J_i(t)$ we
conclude that $S_N(t)J_i(t)=J_i'(t)$ for all $i$.\\

\noindent Now consider the maps
$$J(t):\RR^{n-1}\rightarrow T_{\g(t)}\tau(t),\;\;\;\;\; e_i \mapsto J_i(t)$$ and
$$J'(t):\RR^{n-1}\rightarrow T_{\g(t)}\tau(t),\;\;\;\;\; e_i \mapsto J'_i(t),$$
where $\{e_i\}$ is the standard basis of $\RR^{n-1}$.
As long as the $J_i(t)$'s are linearly independent, we clearly have
$$S_N(t)=  J'(t) \cdot J(t)^{-1}.$$
Propagating the $E_i$'s along $\gamma$ by parallel translation we
obtain an orthonormal basis $\{E_i(t)\}$ of $T_{\g(t)}\tau(t)$.
Furthermore, $\{E_1(t),\dots,E_k(t)\}$ together with $\gd(t)$ span
the vertical space at $\g(t)$ and
$\{E_{k+1}(t),\dots,E_{n-1}(t)\}$ span the horizontal space there.
We will represent the maps $J(t)$,$J'(t)$ and $S_N(t)$ by matrices
w.r.t. the bases $\{e_i\}$ for $\RR^{n-1}$ and $\{E_i(t)\}$  for
$T_{\g(t)}\tau(t)$.\\

\noindent Now we use Jacobi field estimates as in  [BK 6.3.8iii] to determine the operator norm of
$J(t)$, or rather of the endomorphisms $J(t)_{VV}$, $J(t)_{HV}$,
$J(t)_{VH}$ and $J(t)_{HH}$ that $J(t)$ induces on horizontal and vertical
subspaces. \footnote{To be more precise: $J(t)_{HV}: \RR^k \times \{0\}
\rightarrow Hor(t)$ is given by restricting $J(t)$ and then composing with the orthogonal projection
onto the horizontal space at $\g(t)$.} This will allow us to obtain corresponding estimates for $J^{-1}(t)$ and $J'(t)$, and
therefore for $S_N(t)$.\\
\noindent
For all $i$ let us define the vector fields $A_i(t)=\para
(J_i(0)+t\cdot J_i'(0))$, where $\para$ denotes parallel translation along $\gamma$.
The map $\RR^{n-1}\rightarrow T_{\g(t)}\tau(t),\; e_i \mapsto A_i(t),$ in matrix form reads
 $$A(t)=\begin{bmatrix} tI&0\\
 0 & I+tB\end{bmatrix}.$$
 For $i\le k$ we have  $J_i(0)=0$ and  $\{J_i'(0)\}$ is an orthonormal set.
 If $(c_1,\cdots,
c_k,0,\cdots,0)$ is a unit vector in $\RR^{n-1}$, we have $|(\sum
c_i J_i)'(0)|=1$, so applying [BK, 6.3.8iii] we obtain $|\sum c_i
(J_i(t)-A_i(t))| \le  \sinh(t)-t$.\\ Similarly, for $i \ge k+1$,
the set $\{J_i(0)\}$ is an orthonormal set and
$J'_i(0)=B(J_i(0))$. Again, if
$(0,\cdots,0,c_{k+1},\cdots,c_{n-1})$ is a unit vector in
$\RR^{n-1}$, since $|(\sum c_i J_i)'(0)|=|B(\sum c_iJ_i(0))|\le
|B|$, we have
 $|\sum c_i(J_i(t)-A_i(t))| \le  \cosh(t)-1+|B|(\sinh(t)-t)$.
Therefore we have $$J(t)-A(t)=:F_1(t)\le\begin{pmatrix}
\sinh(t)-t\; &\; \cosh(t)-1+|B|(\sinh(t)-t) \\
\sinh(t)-t\; &\; \cosh(t)-1+|B|(\sinh(t)-t)
\end{pmatrix} \le
\begin{pmatrix} \frac{1}{5}t^3 & \frac{3}{4} t^2 \\ \frac{1}{5}t^3 & \frac{3}{4} t^2
\end{pmatrix}
 .$$\\
 \noindent
 Now we want to estimate $tJ^{-1}(t)$. Notice that, suppressing the $t$-dependence in the notation,
 $J=A\cdot[I+A^{-1} F_1]$, so that  $$tJ^{-1}=[I+A^{-1}
 F_1]^{-1}\cdot tA^{-1}.$$
Clearly $A$ is invertible and $$t A^{-1}=\begin{bmatrix} I & 0\\ 0 &
t\cdot(I+tB)^{-1}
\end{bmatrix}\le \begin{pmatrix} 1& 0\\0&2t
\end{pmatrix}$$
since we assume $t\le \frac{1}{2|B|}$. We have
$$A^{-1}
F_1\le \begin{pmatrix}\frac{1}{5}t^2 & \frac{3}{4}t \\
\frac{2}{5}t^3 & \frac{3}{2}t^2
\end{pmatrix}.$$
 Clearly \footnote {If $\begin{bmatrix}  A&B \\C & D
\end{bmatrix}\le
\begin{pmatrix}    a&b \\c & d
\end{pmatrix}$ then the full operator norm of the matrix is bounded by
$\sqrt{\max \{ \begin{vmatrix}   a \\c
\end{vmatrix} ,
\begin{vmatrix}   b \\d   \end{vmatrix} \} +ab+cd }
\le
 \sqrt{2} \max\{ \begin{vmatrix}  a \\c   \end{vmatrix} ,
  \begin{vmatrix}  b \\ d    \end{vmatrix} \}$.}
this is less than $ \sqrt{2}\frac{3}{4}t\sqrt{1+4t^2}\le
\frac{3}{2}t<1$
 since $t\le\frac{1}{2}$. Therefore $I+A^{-1}F_1$ is
 invertible and  $[I+A^{-1}F_1]^{-1}=\sum_{j=0}^{\infty}[-A^{-1}F_1]^j.$
Using the above estimate for $A^{-1} F_1$ we have
$$[A^{-1}F_1]^2\le \begin{pmatrix} \frac{1}{2}t^4 & \frac{3}{2}t^3\\
t^5 &  3t^4
\end{pmatrix}.$$ Using the coarse estimate $A^{-1}
F_1\le \frac{3}{2}t$ and $t\le \frac{1}{2}$ we have
$\sum_{j=3}^{\infty}[-A^{-1}F_1]^j \le 14t^3.$ Putting together these estimates we obtain
$$[I+A^{-1}F_1]^{-1}=I+F_2 \text{ where }F_2\le \begin{pmatrix}
 \frac{15}{2}t^2 & 5t\\
15t^3 &  \frac{19}{2}t^2
\end{pmatrix}.$$

\noindent
To estimate $J'(t)$ we first estimate $|J''(t)-A''(t)|$ and then
integrate. For all $i$ we have
$$|J_i''(t)-A_i''(t)|=|J_i''(t)|\le|J_i(t)|$$
 by the Jacobi equation using the bound on curvature,
and an analogous estimate holds for linear combinations $\sum c_i
J_i(t)$.\\
If $(c_1,\cdots, c_k,0,\cdots,0)$ is a unit vector in $\RR^{n-1}$
we have
$| \sum c_i
J_i(t)|\le \sinh(t)$ by Rauch's theorem.\\
Similarly, if $(0,\cdots,0,c_{k+1},\cdots,c_{n-1})$ is a unit
vector in $\RR^{n-1}$ we have
 $|\sum c_i
J_i(0)|=1$ and  $|\sum c_i J'_i(0)|\le |B|$, so by Berger's
extension of Rauch's theorem (see Lemma 2.7.9 in [Kl]) we have
$|\sum c_i J_i(t)|\le
 \cosh(t)+|B|\sinh(t)$.\\
\noindent In both cases integration delivers \begin{eqnarray*} \arrowvert \sum c_i(J_i'(t)-A'_i(t)) \arrowvert&\le&
 \int_0^t \arrowvert \sum c_i(J_i''(t)-A''_i(t))\arrowvert \text{dt}\\
&\le&
\begin{cases}
\cosh(t)-1\le \frac{3}{4}t^2 &\text{    if }i\le k\\
 \sinh(t)+|B|(\cosh(t)-1)\le \frac{3}{2}t &\text{    if }i \ge k+1.
\end{cases}
\end{eqnarray*}
 So
altogether we obtain $$J'(t)-A'(t)=:F_3(t) \text{ where }F_3(t) \le
 \begin{pmatrix} \frac{3}{4}t^2 & \frac{3}{2}t \\ \frac{3}{4}t^2 & \frac{3}{2}t
\end{pmatrix}.$$\\
Now finally we can estimate
\begin{eqnarray*}
 t S_N(t) &=&tJ'J^{-1}\\
&=&(A'+F_3)\cdot(I+F_2)\cdot tA^{-1}\\
&\le&  \left\{
\begin{bmatrix} I&0 \\0 &B
\end{bmatrix}+
\begin{pmatrix} \frac{3}{4}t^2 &\frac{3}{2}t \\ \frac{3}{4}t^2 &\frac{3}{2}t
\end{pmatrix}+
\begin{pmatrix}  \frac{15}{2}t^2& 5t \\
15|B|t^3
&\frac{19}{2}|B|t^2
\end{pmatrix}+
\begin{pmatrix} 30t^4 & 18t^3\\ 30t^4& 18t^3
\end{pmatrix} \right\}\cdot tA^{-1}\\
&\le&
 \begin{bmatrix}  I&0 \\0 &tB
\end{bmatrix}+
\begin{pmatrix}  0&0 \\0 & 2|B|^2 t^2
\end{pmatrix}+
\begin{pmatrix} \frac{3}{4}t^2 & {3}t^2\\ \frac{3}{4}t^2 & {3}t^2
\end{pmatrix}+
\begin{pmatrix}  \frac{15}{2}t^2& 10t^2 \\15|B| t^3
&19|B|t^3
\end{pmatrix}+
\begin{pmatrix} 30t^4 & 36t^4\\ 30t^4& 36t^4
\end{pmatrix} \\
\end{eqnarray*}

\noindent Here we used
 $$tA^{-1}\le \begin{bmatrix}  I& 0\\0 &tI
\end{bmatrix} + \begin{pmatrix}  0& 0\\0
& 2|B|t^2
\end{pmatrix}
\le  \begin{pmatrix}  1& 0\\0
& 2t
\end{pmatrix}$$ in the last inequality.
Using our bounds on $t$ and the fact that $S_N(t)$ is a symmetric operator this gives the claimed
estimate.
 \Kasten

\noindent Returning to the case when $(M,N)$ is gentle pair, so
that $|B|\le \frac{3}{2}$, we obtain our improvement of Corollary
3.13 in [We]. Now we can achieve an upper bound proportional to
$t^2$, versus the bound proportional to $\sqrt{t}$ of Corollary
3.13 in [We].

\begin{Cor}\label{313}
Let $M$ be a Riemannian manifold, $N$ a submanifold so that
$(M,N)$ form a gentle pair.
If $v$ is a horizontal vector and $w$ a vertical vector at some
point of distance $t\le \frac{1}{3}$ from $N$, then $|\langle
H_N(v), w \rangle| \le 16 t^2|v||w|$.
\end{Cor}

\subsection{Improvement of Weinstein's averaging theorem}\label{s22}

Now we use Corollary \ref{313} to replace some estimates in [We]
that were originally derived using Corollary 3.13 there. We will
improve only estimates contained in Lemma 4.7 and Lemma 4.8 of
[We], where the author considers the covariant derivative of a
certain vector field ${\cal V}$ on $M$ in directions which are
almost vertical or almost horizontal \footnote{See our Section
\ref{s3} or Section 3.2 in [We] for the definitions of almost
horizontal and almost vertical bundle.}
 with respect to a
fixed submanifold $N_e$. (The zero set of ${\cal V}$ is the average $N$ of the family $\{N_g\}$). As in [We] all estimates will hold for
$\epsilon< \frac{1}{20000}$,
 and we set
 $t=100\epsilon$.\\

\noindent We will replace the constant ``$\frac{89}{200}$'' in Lemma 4.7 by
``$\frac{4}{5}$'' as follows:
\begin{Lem}\label{Lemk2}
For any almost vertical vector $v$ at any point of $N$,
$$\langle D_v{\cal V},v \rangle \ge \frac{4}{5}||v||^2.$$
\end{Lem}
\textit{Proof: } By Theorem 3 (applied to the gentle pair $(M,N_g)$) for the operator norm of $H_g$ we
have $1-16t^2\le |H_g|$, so that one obtains
$H_g(\PP_{\Gamma_g}v,\PP_{\Gamma_g}v)>\frac{19}{20}$ in the proof
of Lemma 4.7 in [We]. Similarly, Theorem 3 together with
footnote 11 imply that $|H_g| < 1.01$. Using
these estimates in the proof of Lemma 4.7 in [We] gives the claim.
\Kasten

\noindent Similarly, we will replace the term ``$60 \sqrt{\epsilon}$'' in
Lemma 4.8 by ``$1950 \epsilon$''.
\begin{Lem}\label{Lemk3}
For any almost horizontal vector $v$ at any point of $N$
$$|| D{\cal V}(v) || \ge 1950 \epsilon ||v||.$$
\end{Lem}
\textit{Proof: } By Corollary \ref{313} we can replace ``$3
\sqrt{t}$'' by ``$16 t^2$'' in the proof of Lemma 4.8 in [We] and
we can use 1.01 instead of 1.32 as an upper bound for $|H_g|$.
Furthermore, we replace the constant 1000 coming from Lemma 4.3 in
[We] by 525.\footnote{Lemma 4.3 of [We] quotes incorrectly
Proposition A.8 from its own appendix.} This delivers the improved
estimate $H_g(v,\PP_{\bar{\Gamma}}w)\le 850 \epsilon||v||\cdot
||w||$ and simple arithmetics concludes the proof. \Kasten

\noindent From these two lemmas it follows that the operator from
$(aVert^e)^{\perp}$ to $aVert^e$ whose graph is $T_xN$  has norm at
most $\frac{5}{4}\cdot 1950\epsilon$. Following to the end the
proof of Theorem 2.3 in [We] allows us to replace the bound
``$136\sqrt{\epsilon}$'' by a bound linear in $\epsilon$, so we
obtain the following improved statement:

\begin{Thm}\label{Thmk}
 Let $M$ be a \riem \mfd and
$\{N_g\}$ a family of submanifolds of $M$ parametrized in a
measurable way by elements of a probability space $G$, such that
all the pairs $(M,N_g)$ are gentle. If
$d_1(N_g,N_h)<\epsilon<\frac{1}{20000}$ for all $g$ and $h$ in
$G$, there is a well defined \emph{\textbf{center of mass}}
submanifold $N$ with $d_1(N_g,N)<2500 \epsilon$ for all $g$ in
$G$. The center of mass construction is equivariant with respect
to isometries of $M$ and measure preserving
automorphisms of $G$.\\
\end{Thm}

\section{Estimates on the map $\vg$}\label{s3}

\noindent In sections \ref{s3}-\ref{s7} we will prove the Main Theorem. The reader is referred to Section \ref{s13} for an outline of the proof and some of the notation introduced there. We will present Part I of the proof in sections
\ref{s3}-\ref{s6} and Part II in Section \ref{s7}.\\

 Fix $g\in G$ and let $p$ be a
point in the tubular \nbhd of $N_g$ and $X \in T_pM$. The aim of
this section
is to estimate the difference between $\vg_*X$ and $\para X$. This will be achieved in Proposition 3.4.\\
Here we denote by $\para X$ the following parallel translation of
$X$, where $\png$ is the projection onto $N_g$ along the normal
slices. First we parallel translate $X$ along the shortest
geodesic from $p$ to $\png(p)$, then along the shortest geodesic
from $\png(p)\in N_g$ to its image under $\vg$, and finally along
the shortest geodesic to $\vg(p)$. We view ``$\para$'' as a
canonical way to associate $X\in T_pM$ to a vector in
$T_{\vg(p)}M$.\\[0.4cm]
Before we begin proving our estimates, following
section 2.1 of [We] we introduce two subbundles of
$TM|_{\tnng}$ and their orthogonal complements.\\
The \textit{vertical bundle} $\vertg$ has fiber at $p$
given by the parallel translation of
$\nu_{\png(p)}N_g$ along the shortest geodesic from
$\png(p)$ to $p$.\\
The \textit{almost vertical bundle} $\avg$ has fiber
at $p$ given by the tangent space at $p$ of the normal
slice to $N_g$ through $\png(p)$.\\
The  \textit{horizontal bundle} $\horg$ and the
\textit{almost horizontal bundle} $\ahg$ are given by
their orthogonal complements.\\[0.4cm]
\rem Notice that $\avg$ is the kernel of $(\png)_*$,
and that according to Proposition 3.7 in [We] we have
$d(\vertg_p,\avg_p)<\frac{1}{4}d(p,N_g)^2$ for any $p$
in $\tnng$, and similarly for $\horg$ and $\ahg$.\\
Since $\frac{1}{4}d(p,N_g)^2<\frac{\pi}{2}$, $\vertg$ and $\ahg$
are always transversal (and clearly the same holds for $\horg$ and
$\avg$). As seen in Section \ref{s13} $\vertg$ and $TN$ are
transversal along $N$, and $\avg$ and $TN$ are also transversal
since $N$ corresponds to a section of $\nu N_g$ and
$\avg=\text{Ker}(\pi_{N_g})_*$.\\[0.4cm]
Now we are ready to give our estimates on the map $\vg$. Recall
from the Introduction that for any point $q$ of the tubular \nbhd
of $N_g$ we denote by $\gamma_q$ the geodesic from $\png(q)\in
N_g$ to $q$. Until the end of this section all geodesics will be
parametrized to arc-length.\\In sections \ref{s3} to \ref{s6} all
estimates will hold for $\e < \frac{1}{20000}$.\\

\subsection{Case 1:  $p$ is a point of $N_g$}\label{s31}

\begin{Prop}\label{Prop21}  If $p\in N_g$ and $X\in T_pN_g$ is a
unit vector, then $$\left|\vg_*(X) - \para X\right| \le 3200 \e.$$
\end{Prop} \rem Notice that if $X$ is a vector normal to $N_g$ by
definition of $\vg$ and $\para$ we have $\vg_*(X)=\para X$.
Therefore in this subsection we will assume that $X$ is tangent to
$N_g$.
\\[0.4cm]
Also, we will denote by $A$ the second fundamental form
\footnote{In Section \ref{s2} we adopted the sign convention of
[We] which differs from this.}
 $N_g$, i.e. $A_{\xi}v:= -(\nabla_v \xi)^T$ for tangent
vectors $v$ of $N_g$ and normal \vfs $\xi$, where $(\cdot)^T$
denotes projecting to the component tangent to $N_g$ and $\nabla$
is the Levi-Civita connection on $M$. Since $(M,N_g)$ is a gentle
pair, the norm of $A$ is bounded by $\frac{3}{2}$, as shown
in [We, Cor. 3.2].\\[0.4cm]
Now let $p\in N_g$, $X \in T_p N_g$ a unit vector, and
$q:=\vg(p)$. We will denote by $E$ the distance $d(p, \vg(p))< 100
\e$ (see end of Section 4 in [We]).\\ We will show that at
$q$
$$\ \para X \approx J(E) \approx H \approx \vg_*(X) $$ where the
Jacobi field $J$ and the horizontal vector $H$ will be specified
below.

\begin{Lem}\label{Lem21} Let $J$ be the \jf along the geodesic
$\gamma_{q}$ such that $J(0)=X$ and
$J'(0)=-A_{\dot{\gamma}_{q}(0)}X$. Then $$\left| J(E)-\para X
\right|\le \frac{3}{2}(e^E-1).$$
\end{Lem}
\pf
This is an immediate consequence of [BK,6.3.8.iii]
\footnote{[BK,6.3.8] assumes that $J(0)$ and $J'(0)$ be linearly
dependent. However statement iii) holds without this assumption,
as one can always decompose $J$ as $J=J_1+J_2$, where $J_1$ and
$J_2$ are Jacobi fields  such that $J_1(0)=J(0),J_1'(0)=0$ and
$J_1(0)=0,J_1'(0)=J'(0)$ respectively.
  Furthermore we make use of $|J|'(0)\le
|J'(0)|$.}
 which will be
used later again and which under the curvature assumptions $|K|\le
1$ states the following: if $J$ is any Jacobi field along a
unit-speed geodesic, then we have
$$\left| J(t)- _{t}^{0}\para \left(J(0)+t\cdot J'(0) \right)
\right|\le
\left|J(0)\right|(\cosh(t)-1)+\left|J'(0)\right|(\sinh(t)-t ),$$
where $_{t}^{0}\para $ denotes parallel translation to the
starting point of the geodesic.
 Using $| A_{\dot{\gamma}_q(0)} X| \le \frac{3}{2}$ by [We,
Cor 3.2] the above estimates gives $\left| J(E)-\para X \right|\le
(\cosh(E)-1)+ \frac{3}{2}\sinh(E)$. Alternatively, this Lemma can
be proven using the methods of [We, Prop 3.7].
\Kasten\\[0.2cm]

\noindent Before proceeding we need a lemma about projections:

\begin{Lem}\label{Lem22} If $Y\in T_qM$ is a vertical unit vector, write
$Y=Y_{av} +Y_h$ for the splitting into its almost vertical and
horizontal components. Then
$$\left|Y_h\right|\le \tan\left(\frac{E^2}{4}\right)\;\;\;\text{
and }\;\;\;\left|Y_{av}\right|\le \frac{1}{\cos(\frac{E^2}{4})}.
$$

\end{Lem} \pf By [We, Prop 3.7]
we have $d(\vertg_q,\avg_q) \le \frac{E^2}{4} < \frac{\pi}{2}$, so
the subspace $\avg_q$ of $T_qM$ is the graph of a linear map
$\phi: \vertg_q \rightarrow \horg_q$. So $Y_{av}=Y+\phi(Y)$ and $
Y_h=-\phi(Y)$. Since the angle enclosed by $Y$ and $Y_{av}$ is at
most $d(\vertg_q,\avg_q) \le \frac{E^2}{4}$, one obtains $|Y| \ge
\cos(\frac{E^2}{4})|Y_{av}|$ which gives the second estimate of
the Lemma. From this, using $|Y_h|^2=|Y_{av}|^2 -|Y|^2$ we obtain
the first estimate. \Kasten\\[0.2cm]

\begin{Lem}\label{Lem23} If $H$ is the unique horizontal vector at
$q$ \st ${\png}_*(H)=X$, then $$\left|  J(E) -H \right| \le
\frac{3}{2}(e^E-1)\frac{1}{\cos (\frac{E^2}{4})}.$$\end{Lem} \pf
Let $J$ be the \jf of Lemma \ref{Lem21}. Write $J(E)=W +Y$ for the
splitting into horizontal and vertical components. Then, using the
notation of Lemma \ref{Lem22}, we have $J(E)_h=W+Y_h$ and
$J(E)_{av}=Y_{av}$. Notice that the Jacobi field $J$ arises from a
variation of geodesics \orth to $N_g$ (see the Remark in Section
\ref{s32}), so $(\png)_*J(E)=X= (\png)_*H$. Using
$\avg=\ker(\png)_*$ it follows that $H=J(E)_h$. So
$$|J(E)-H|=|Y_{av}|\le |Y| \frac{1}{\cos (\frac{E^2}{4})}
 \le
\frac{3}{2}(e^E-1)\cdot \frac{1}{\cos (\frac{E^2}{4})}$$ where we
used Lemma \ref{Lem22} and $\left|Y\right| \le \left|J(E) -
\para X \right|$ together with Lemma \ref{Lem21}.
\Kasten\\[0.2cm]
\noindent Now we will compare $H$ to $\vg_*(X)$
and finish our proof.\\[0.2cm]
\textit{Proof of Proposition  \ref{Prop21}: }
 We have $$ \left|  \para X - \vg_*(X)
\right| \le \left| \para X -J(E) \right| + \left|J(E) -
H\right|+\left|H-\vg_*(X)\right|.$$ The first and second terms are
bounded by the estimates of Lemmata \ref{Lem21} and \ref{Lem23}.
For the third term we proceed analogously to Lemma \ref{Lem23}:
since $\vg_*(X)$ and $H$ are both mapped to $X$ via $\png$, one
has $(\vg_*(X))_{av}=\vg_*(X)-H$. As earlier, if $\vg_*(X)=
\tilde{W} +\tilde{Y}$ is the splitting into
 horizontal and vertical
components, we have $(\vg_*(X))_{av}=\tilde{Y}_{av}$.
 Therefore
$$|\vg_*(X)-H|=|\tilde{Y}_{av}| \le |\tilde{Y}| \frac{1}{\cos(\frac{E^2}{4})} \le |\vg_*(X)|
\frac{\sin (2500 \e)}{\cos(\frac{E^2}{4})}.$$ Here we also used
Lemma \ref{Lem22} and the fact that the angle enclosed by
$\vg_*(X)$ and its orthogonal projection onto $Hor_g^q$ is at most
$d(\horg_q, T_qN) \le 2500 \e$ by Theorem 4.
Altogether we have $$\left|  \para X - \vg_*(X) \right| \le
\frac{3}{2}(e^E-1)\left[1+\frac{1}{\cos(\frac{E^2}{4})}
\right]+|\vg_*(X)|\frac{\sin (2500 \e)}{\cos(\frac{E^2}{4})}.$$
 Using
this inequality we can bound $|\vg_*(X)|$ from above in terms of
$E$ and $\e$. Substituting into the right hand side of the above
inequality we obtain a function of $\e$ (recall that $E=100
\e$)which is increasing and bounded above by $3200 \e$. \Kasten
%
\subsection{Case 2: $p$ is a point of $\partial
\tnngl$ and $X\in T_pM$ is almost vertical}\label{s32}

In this subsection we require $L<1$, as in the
definition of gentle pair.\\[0.4cm]
\rem Jacobi-fields $\bar{J}$ along $\gamma_p$ (the geodesic from
$\png(p)$ to $p$) with $\bar{J}(0)$ tangent to $N_g$ and $
A_{\dot{\gamma}_p(0)}\bar{J}(0)+\bar{J}'(0)$ normal to $N_g$ are
called \textit{$N_g$ Jacobi-fields}. They clearly form a \vs of
dimension equal to $\dim(M)$ and they are exactly the \jfs that
arise from variations of $\gamma_p$ by geodesics that start on
$N_g$ and are normal to $N_g$ there.\\Since $(M,N_g)$ is a gentle
pair, there are no focal points of $\png(p)$ along $\gamma_p$, so
the map $$\{N_g \text{ \jfs along } \gamma_p \} \rightarrow
T_pM,\: \bar{J} \mapsto \bar{J}(L)$$ is an isomorphism. The $N_g$
Jacobi fields that map to $aVert_p^g$ are exactly those with the
property $J(0)=0$, $J'(0)\in \nu_{\png(p)}N_g$. Indeed such a \vf
is the variational \vf of the variation
$$f_s(t)=\exp_{\png(p)}t[\dot{\gamma_p}(0)+sJ'(0)],$$
so $J(L)$ will be tangent to the normal slice of $N_g$ at
$\png(p)$. From dimension considerations it follows that the $N_g$
\jfs that satisfy $J(0)\in T_{\png(p)}N_g$ and
$A_{\dot{\gamma}_p(0)}J(0) + J'(0)=0$  - which are called \textit{
strong $N_g$ \jfs} - map to a subspace of $T_pM$ which is a
complement of $aVert_p^g$. As pointed out in [Wa, p. 354], these
two subspaces are in general not
orthogonal.\\
\begin{Prop}\label{Prop22} If $p\in \partial \tnngl$ and $X\in T_pM$
is an almost vertical unit vector, then
$$\left|\vg_*(X)- \para X\right| \le
2\frac{\sinh(L)-L}{\sin(L)}.$$ \end{Prop} We begin proving
\begin{Lem}\label{Lem24} Let $J$ be a \jf along $\gamma_p$ \st
$J(0)=0$ and $J'(0)\in \nu_{\png(p)} N_g$, normalized \st
$|J(L)|=1$. Then $$\left|J(L)-L \cdot _{\gamma_p}\parap
J'(0)\right|\le \frac{\sinh(L)-L}{\sin(L)}.$$ \end{Lem} \pf Again
[BK, 6.3.8iii] shows that $\left|J(L)-L \para J'(0)\right|\le
|J'(0)|(\sinh(L)-L).$ Using the upper curvature bound $K\le 1$ and
Rauch's theorem we obtain $|J'(0)| \le \frac{1}{\sin(L)}$ and we
are done.  \Kasten\\[0.2cm]\\
We saw in the remark above that $X$ is equal to $J(L)$ for a \jf
$J$ as in Lemma \ref{Lem24}, and that $J$ comes from a variation
$f_s(t)=\exp_{\png(p)}t[\dot{\gamma_p}(0)+sJ'(0)]$. So $\vg_*(X)$
comes from the variation $$\vg \left(f_s(t)
\right)=\exp_{\vg(\png(p))} t\left[\para \dot{\gamma_p}(0)+s\para
J'(0)\right]$$ along the geodesic $\vg(\gamma_p(t))$. More
precisely, if we denote by $\tilde{J}(t)$ the \jf that arises from
the above variation, we will have $\vg_*(X)=\tilde{J}(L)$. Notice
that $\tilde{J}(0)=0$
and $\tilde{J}'(0)= \para J'(0)$.\\
\begin{Lem}\label{Lem25} $$\left| \tj(L) -L \cdot  _{\vg \circ
\gamma_p}\parap \tj '(0) \right| \le \frac{\sinh(L)-L}{\sin(L)}.$$
\end{Lem} \pf Exactly as for Lemma \ref{Lem24} since $\tilde{J}(0)=0$ and
$|\tilde{J}'(0)|=|J'(0)|.$
 \Kasten\\[0.2cm]
\textit{Proof of Proposition \ref{Prop22}: } We have $X \approx
LJ'(0) = L \tilde{J}'(0) \approx \vg_*(X)$. Here we identify
tangent spaces to $M$ parallel translating along $\gamma_p$, along
the geodesic $\gamma_{\vg(\png(p))}$ from $\png(p)$ to its
$\vg$-image and along $\vg \circ\gamma_p$ respectively. Notice
that these three geodesics are exactly those used in the
definition of
``$\para$''.\\
The estimates for the two relations ``$\approx$'' are in Lemma
\ref{Lem24} and Lemma \ref{Lem25} respectively (recall $X=J(L)$
and $\vg_*(X)=\tilde{J}(L)$), and the equality holds because
$\tilde{J}'(0) =\para J'(0).$\Kasten

\subsection{Case 3:  $p$ is a point of $ \partial
\tnngl$ and $X=J(L)$ for some  strong $N_g$ \jf $J$ along
$\gamma_p$}\label{s33}

From now on we have to assume $L < 0.08$.
\begin{Prop}\label{Prop23} If $p\in \partial \tnngl$ and $X$ is a
unit vector equal to $J(L)$ for some strong $N_g$ \jf $ J$ along
$\gamma_p$, then $$\left|\vg_*(X)- \para X \right| \le
\frac{18}{5}L+3700 \e.$$ \end{Prop} We proceed analogously to Case
2.
\begin{Lem}\label{Lem26} For a \vf $J$ as in the above proposition we have $$|J(L) -
_{\gamma_p} \parap J(0)|\le \frac{\frac {3}{2}(e^L-1)}{1- \frac
{3}{2}(e^L-1)}\le \frac{9}{5}L.$$ Furthermore we have $|J(0)| \le
\frac{1}{1- \frac {3}{2}(e^L-1) }.$ \end{Lem} \pf By Lemma
\ref{Lem21} we have $|J(L) - _{\gamma_p} \parap J(0)| \le
\frac{3}{2} (e^L -1)|J(0)|$, from which we obtain the estimate for
$|J(0)|$ and then the first
estimate of the lemma. \Kasten\\[0.2cm]
$J$ comes from a variation
$f_s(t)=\exp_{\sigma(s)}tv(s)$ for some curve $\sigma$
in $N_g$ with $\dsigma(0)=J(0)$ and some normal \vf
$v$ along $\sigma$. We denote by $\tj$ the \jf along
the geodesic $\vg(\gamma_p(t))$ arising from the
variation $$\tilde{
f_s}(t)=\vg
\left(f_s(t)\right)=\exp_{\tilde{\sigma}(s)} \left(t
\para v(s) \right),$$ where $\tilde{\sigma}=\vg
\circ \sigma$ is the lift of $\sigma$ to $N$. Then we
have $\tj(L)=\vg_*(X)$. Notice that here $\para
v(s)$ is just the parallel translation of $v(s)$ along
$\gamma_{\tilde{\sigma}(s)}=:\gamma_s$.\\
\begin{Lem}\label{Lem}$$\left|\tj(L)- _{\vg \circ \gamma_p}
\parap \tj(0) \right |\le \frac{9}{5}L$$ \end{Lem}
\pf
Using [BK 6.3.8iii] as in Lemma \ref{Lem21} we obtain
 \begin{equation*} \begin{split} \left|\tilde{J}(L) - _{\vg \circ
\gamma_p} \parap \tilde{J}(0)\right|
\le& |\tilde{J}'(0)| ( \sinh(L)) + |\tilde{J}(0)|(\cosh(L)-1 ),
\;\;\;\;\;\;(*)
\end{split} \end{equation*}
so that we just have to estimate the norms of
$\tilde{J}(0)$ and $\tilde{J}'(0)$.\\
Since $\tilde{J}(0)= \vg_*J(0)$, applying Proposition \ref{Prop21}
gives $|\tilde{J}(0) - _{\gamma_0}\parap J(0)| \le 3200 \e
|J(0)|$. Using the bound for $|J(0)|$ given in Lemma \ref{Lem26}
we obtain $|\tilde{J}(0)|\le
\frac{1+3200\e}{1- \frac{3}{2} (e^L -1)}$.\\
To estimate $\tilde{J}'(0)$ notice that in the
expression for $f_s(t)$ we can choose $v(s)=
_{\sigma_s}\parap[ \dot{\gamma}_0(0)+sJ'(0)]$, where
$_{\sigma_s}\parap$ denotes parallel translation from
$\sigma(0)$ to $\sigma(s)$ along $\sigma$. So $$ \para
v(s)= _{\gamma_s}\parap _{\sigma_s}\parap\left[
\dot{\gamma}_0(0)+sJ'(0) \right],$$ and $$
\tilde{J}'(0)= \frac{\nabla}{ds}\Big|_0 \big(\para
v(s) \big)=\frac{\nabla}{ds}\Big|_0\,
_{\gamma_s}{\parap} _{\sigma_s}{\para}
\dot{\gamma}_0(0) + _{\gamma_0}\parap J'(0)$$ where we
used the Leibniz rule for covariant derivatives to
obtain the second equality.\\
To estimate the first term note that the difference between the
identity and the holonomy around a loop in a Riemannian manifold
is bounded in the operator norm by the surface spanned by the loop
times the a bound for the curvature (see [BK, 6.2.1]). Therefore
we write $ _{\gamma_s}{\parap} _{\sigma_s}{\parap}
\dot{\gamma}_0(0)$ as $ _{\tilde{\sigma}_s}
\parap _{\gamma_0}
\parap
 \dot{\gamma}_0(0) + \varepsilon(s)$
where $\varepsilon(s)$ is a \vf along $\tilde{\sigma}(s)$ with
norm bounded by the area of the polygon spanned by $\sigma(0),
\sigma(s), \tilde{\sigma}(s)$ and $\tilde{\sigma}(0)$. Assuming
that $\sigma$ has constant speed $|J(0)|$ we can estimate
$d(\sigma(0), \sigma(s)) \le s|J(0)|$ and using Proposition
\ref{Prop21} to estimate $|\dot{\tilde{\sigma}}(s)| =|\vg_*
\dot{\sigma}(s)|$ we obtain  $d(\tsigma(0), \tsigma(s)) \le
s(1+3200 \e)|J(0)|$. Using $d(\tsigma(s), \sigma(s)) \le 100 \e$
and Lemma \ref{Lem26} we can bound the area of the polygon safely
by $\frac{100\e s(2+3200\e)}{1- \frac{3}{2} (e^L -1)}.$ So we
obtain
$$\left| {\frac{\nabla}{ds}\Big|_0}
\,_{\gamma_s}{\parap} _{\sigma_s}\parap \dot{\gamma}_0(0)\right|
=\left|\frac{\nabla}{ds}\Big|_0 \varepsilon(s)\right| \le
\frac{100\e(2+3200\e)}{1- \frac{3}{2} (e^L -1)}.$$ To bound $
_{\gamma_0}\parap J'(0)$ notice that $|J'(0)| \le \frac{3}{2}
|J(0)|$ using the fact that $J$ is a strong \jf and [We, Cor.
3.2], so $|J'(0)|
\le \frac{3}{2} \frac{1}{1- \frac{3}{2} (e^L -1)}.$\\
Substituting our estimates for  $|\tilde{J}(0)|$ and
$|\tilde{J}'(0)|$ in $(*)$ we obtain a function which,
for $\e < \frac{1}{20000}$ and $L < 0.08$, is bounded above by
$\frac{9}{5}L$.\Kasten\\[0.2cm]
The vectors $\tilde{J}(0)$ and $\para  J(0)$
generally are not equal, so we need one more estimate
that has no counterpart in Case 2:
\begin{Lem}\label{Lem28}$$\left|\tj(0) - \para J(0) \right| \le
\frac{3200 \e}{1- \frac {3}{2}(e^L-1)}\le 3700 \e.$$
\end{Lem}
\pf Since $\tilde{J}(0)= \vg_*J(0)$, Proposition  \ref{Prop21}
gives
$$|\tilde{J}(0) - \para J(0)| \le 3200 \e|J(0)| \le
\frac{3200\e}{1- \frac{3}{2} (e^L -1)}.$$ Since $\frac{1}{1-
\frac{3}{2} (e^L -1)} < 1.15$ when
$L<0.08$ we are done. \Kasten\\[0.2cm]
\textit{Proof of Proposition  \ref{Prop23}: } We have $X \approx
J(0) \approx \tilde{J}(0) \approx \vg_*(X)$ where we identify
tangent spaces by parallel translation along $\gamma_p$,
$\gamma_0$ and $\vg \circ \gamma_p$ respectively. Combining the
last three lemmas and recalling $X=J(L)$, $\vg_*(X)=\tilde{J}(L)$
we finish the proof. \Kasten

\subsection{The general case}\label{s34}

\noindent This proposition summarizes the three cases considered
up to now:
\begin{Prop}\label{Prop24} Assume $\e<\frac{1}{20000}$ and $L<0.08$. Let
$p\in \partial \tnngl$ and $X\in T_pM$ a unit vector. Then $$
\left|\vg_*(X)- \para X \right| \le 4 L + 4100 \e.$$ \end{Prop}

\noindent We will write the unit vector $X$ as $J(L)+K(L)$ where
$J$ and $K$, up to normalization, are Jacobi fields as in the next
lemma. We will need to estimate the norms of $J(L)$ and $K(L)$, so
we begin by estimating the angle they enclose:
\begin{Lem}\label{Lem29} Let $J$ be a  $N_g$ \jf along $\gamma_p$
with $J(0)=0$, $J'(0)$ normal to $N_g$ (as in Case 2) and $K$ a
strong $N_g$ \jf (as in Case 3), normalized \st $J(L)$ and $K(L)$
are unit vectors. Then $$ \left|\la J(L), K(L) \ra \right| \le
\frac{\frac {3}{2}(e^L-1)}{1- \frac {3}{2}(e^L-1)}+ \frac{1}{1-
\frac {3}{2}(e^L-1)}\cdot \frac{\sinh(L)-L}{\sin(L)}\le
\frac{9}{5}L.$$
\end{Lem} \pf Identifying tangent spaces  along $\gamma_p$ by
parallel translation, we have
\begin{equation*} \begin{split} \left| \la J(L),K(L)
\ra \right| =&  \left| \la J(L),K(L) \ra -   \la L J'(0),K(0) \ra
\right| \\ \le&  \left| \la J(L),K(L)-K(0) \ra \right| +  \left|
\la J(L)-LJ'(0),K(0) \ra \right| \\ \le& \left|K(L)-K(0)\right| +
\left|K(0) \right|\cdot \left|J(L) -LJ'(0)\right|, \end{split}
\end{equation*} which can be estimated using Lemma \ref{Lem26} and
Lemma \ref{Lem24}
\Kasten\\[0.2cm]
\begin{Lem}\label{Lem210} Let $X\in T_pM$ be a unit vector \st
$X=J(L)+K(L)$ where $J,K$ are \jfs as in the Lemma \ref{Lem29} (up
to normalization). Then
$$\left| J(L) \right|, \left| K(L)\right| \le \frac{1}{\sqrt{1
-\frac{9}{5}L}}\le 1.1.$$
\end{Lem}
\pf Let $c:= \la
\frac{J(L)}{|J(L)|},\frac{K(L)}{|K(L)|}  \ra$, so $|c|
\le \frac{9}{5}L$. There is an \onb $\{e_1,e_2\}$ of
span$\{ J(L), K(L) \}$ \st $J(L)=|J(L)|e_1$ and
$K(L)=|K(L)|(ce_1 + \sqrt{1 -c^2}e_2)$. An elementary
computation shows that $1=|J(L)+K(L)|^2 \ge
(1-|c|)(|J(L)|^2 + |K(L)|^2)$, from which the lemma
easily follows. \Kasten\\[0.2cm]
\textit{Proof of Proposition  \ref{Prop24}: } The remark at the
beginning of \textit{Case 2 } implies that we can (uniquely) write
$X=J(L)+K(L)$ for $N_g$ Jacobi-fields $J$ and $K$ as in Lemma
\ref{Lem210}. So, using Lemma \ref{Lem210}, Proposition
\ref{Prop22} and Proposition  \ref{Prop23}
\begin{equation*} \begin{split}
\left|\vg_*(X) - \para X\right| \le& \left|\vg_*J(L)
-\para J(L)\right| +  \left|\vg_*K(L) -\para K(L)
\right| \\
\le& 1.1 \left (2\frac{\sinh(L)-L}{\sin(L)} +
\frac{18}{5}L + 3700\e \right)\\
\le&  4L + 4100 \e.
\end{split} \end{equation*} \Kasten

\section{Proposition 4.1 about geodesic triangles in $M$}\label{s4}

Fix $g$ in $G$ and let $\vg$ be the map from a tubular \nbhd of
$N_g$ to one of $N$ defined in the introduction. Our aim in the
next three sections is to show that $\exp_{N_g}(\nu N_g)_{0.05}$
is a tubular
\nbhd of $N_g$ on which $\vg$ is injective.\\
We will begin by giving a lower bound on the length of the edges
of certain geodesic \trs in
$M$.\\[0.4cm]
In this section we take $M$ to be simply any Riemannian \mfd with
the following two properties: \begin{quote}
i) the sectional curvature lies between -1 and 1\\
ii) the injectivity radius at any point is at least 1.
\end{quote}
In our later applications we will work in the \nbhd of
a \subm that forms a gentle pair with $M$, so these
two conditions will be  automatically
satisfied.\\[0.2cm]
Now let us choose points $A,B,C$ in $ M$ and let us assume
$d(C,A)<0.15$ and $d(C,B)<0.5$. Connecting the three points by the
unique shortest geodesics defined on the interval $[0,1]$, we
obtain a \geo \tr $ABC$.
\\We will denote by the symbol $\dot{CB}$ the initial
velocity vector of the \geo from $C$ to $B$, and
similarly for the other edges of the triangle.
\begin{Prop} \label{Prop31} Let $M$ be a Riemannian \mfd and $ABC$ a
\geo \tr as above. Let $P_C,P_A$ be subspaces of $T_CM$  and
$T_AM$ respectively of equal dimensions \st $\dot{CB}\in P_C$ and
$\dot{AB}\in P_A.$  Assume that:
$$\measuredangle(P_A, \dot{AC}) \ge \frac{\pi}{2}-
\delta$$ and  $$\theta:=d(P_A, \,_{CA}\parap P_C) \le {\cal C}
d(A,C)$$ for some constants $\delta, \cal C$. We assume ${\cal
C}\le2$. Then $$d(C,B) \ge \frac{10}{11}\frac{1}{ ({\cal C}+6)
\sqrt{1 + {\tan}^2 \delta}}.$$ \end{Prop} \textbf{Remark 1: }Here
$\:_{CA}\parap$ denotes parallel translation of $P_C$ along the
geodesic from
$C$ to $A$.\\
The \textit{angle} between the subspace $P_A$ and the
vector $\dot{AC}$ is given as follows: for every
nonzero $v\in P_A$ we consider the non-oriented angle
$\measuredangle(v, \dot{AC}) \in [0,\pi].$ Then we
have  $$\measuredangle(P_A, \dot{AC}):=\text{Min
}\left\{ \measuredangle(v, \dot{AC}): v\in P_A \text{
nonzero } \right\}\in \left[0, \frac{\pi}{2}
\right].$$ Notice that $  \measuredangle(P_A,
\dot{AC}) \ge \frac{\pi}{2} - {\delta}$ iff for all
nonzero $v \in P_A$ we have  $\measuredangle(v,
\dot{AC})\in  [\frac{\pi}{2} - {\delta},
\frac{\pi}{2}+ {\delta}].$\\
\textbf{Remark 2: } This proposition generalizes the
following simple statement about triangles in the
plane: if two edges $CB$ and $AB$ form an angle
bounded by the length of the base edge $AC$ times a
constant ${\cal C}$, and if we assume that $CB$ and
$AB$ are nearly perpendicular to $AC$, then the
lengths $|CB|$ and $|AB|$ will be bounded below by a
constant depending on ${\cal C}$ (but not on
$|AC|$).\\
In the general case of Proposition  \ref{Prop31} however we make
assumptions on $d(P_A, \,_{CA}\parap P_C)$ from which we are not
able to obtain easily bounds on the angle
$\measuredangle(\dot{CB},\dot{AB})$ at $B$ (such a bound together
with the law of sines would immediately
imply the statement of Proposition  \ref{Prop31}).\\[0.4cm]
\pf Using the chart $\exp_A$ we can lift $B$ and $C$ to the points
$\tb$ and $\tc$ of $T_AM$. We obtain a triangle $0\tb\tc$, which
differs in one edge from the lift of the triangle $ABC$. Denoting
by  $\tbc$ the endpoint of the vector $(\tb - \tc)$ translated to
the origin, consider the \tr $0\tb \tbc$. Let $P$ be the
closest \pt to $\tbc$ in $P_A$.\\[0.2cm]
\textit{Claim 1: } $$|\tb -P| \le \tan(\delta)|\tbc -
P|.$$\\
Using $\measuredangle(P_A, \dot{AC}) \ge \frac{\pi}{2}
- \delta$ and $\dot{AC}=\tc -0$ we see that the angle
between any vector in $P_A$ and $\tc -0$ lies in the
interval  $ [\frac{\pi}{2} - \delta,  \frac{\pi}{2}+
\delta]$. Since $\tc -0$ and $ \tbc -\tb$ are
parallel, the angle between any vector of $P_A$ and
$\tbc -\tb$ lies in $ [\frac{\pi}{2} - \delta,
\frac{\pi}{2}+ \delta]$. Since $P-\tb\in P_A$ we have
$$\measuredangle\left(P-\tb,\tbc -\tb \right) \in
\left[\frac{\pi}{2} - \delta,  \frac{\pi}{2}+
\delta\right].$$\\
 The triangle $\tb P \tbc$ has a right angle at $P$,
so the angle $ \measuredangle(P- \tbc,\tb- \tbc)$ is
less of equal than $\delta$, and  claim 1
follows.\\[0.2cm]
\textit{Claim 2: }$$ |\tbc -P| \le \sin \left[(6+{\cal
C})\cdot d(C,A)\right]\cdot|\tbc -0|.$$\\
$(\expi_A)_*\dot{CB}$ does not coincide with $\tb -\tc \in
T_{\tc}(T_AM)$. Estimating - see Corollary \ref{CorA2} in Appendix
A - the angle between these two vectors we will be able to
estimate the angle between $\tb-\tc=\tbc -0\in T_AM$ and
$_{CA}\parap \dot{CB}\in T_AM$, i.e. the parallel translation in
$M$ of $\dot{CB}$ along the \geo from $C$ to $A$ (see Corollary
\ref{CorA3}). Our estimate will be
$\measuredangle\left(\,_{CA}\parap \dot{CB},\tbc-0
\right)<4d(A,C).$ \\
Now let $P'$ be the closest point to $_{CA}\parap\, \dot{CB}$ in
$P_A$. Since $P'-0 \in P_A$ and $\dot{CB}\in P_C$, using the
definition of distance between subspaces we get
$\measuredangle(\:_{CA}\parap \dot{CB},P'-0)\le
d(P_A,\;_{CA}\parap P_C )= \theta
\le {\cal C}d(C,A).$\\
Finally we will show - see Corollary \ref{CorA1} - that
$\measuredangle(P-0,P'-0)\le
2d(C,A).$\\
Combining the last three estimates we get $\measuredangle(P-0,
\tbc -0)< (6+{\cal C}) d(C,A)$, which is less than
$\frac{\pi}{2}$.\\ Claim 2 follows since $0P \tbc$ is a right \tr
at $P$.\\

\centerline{\epsfig{file=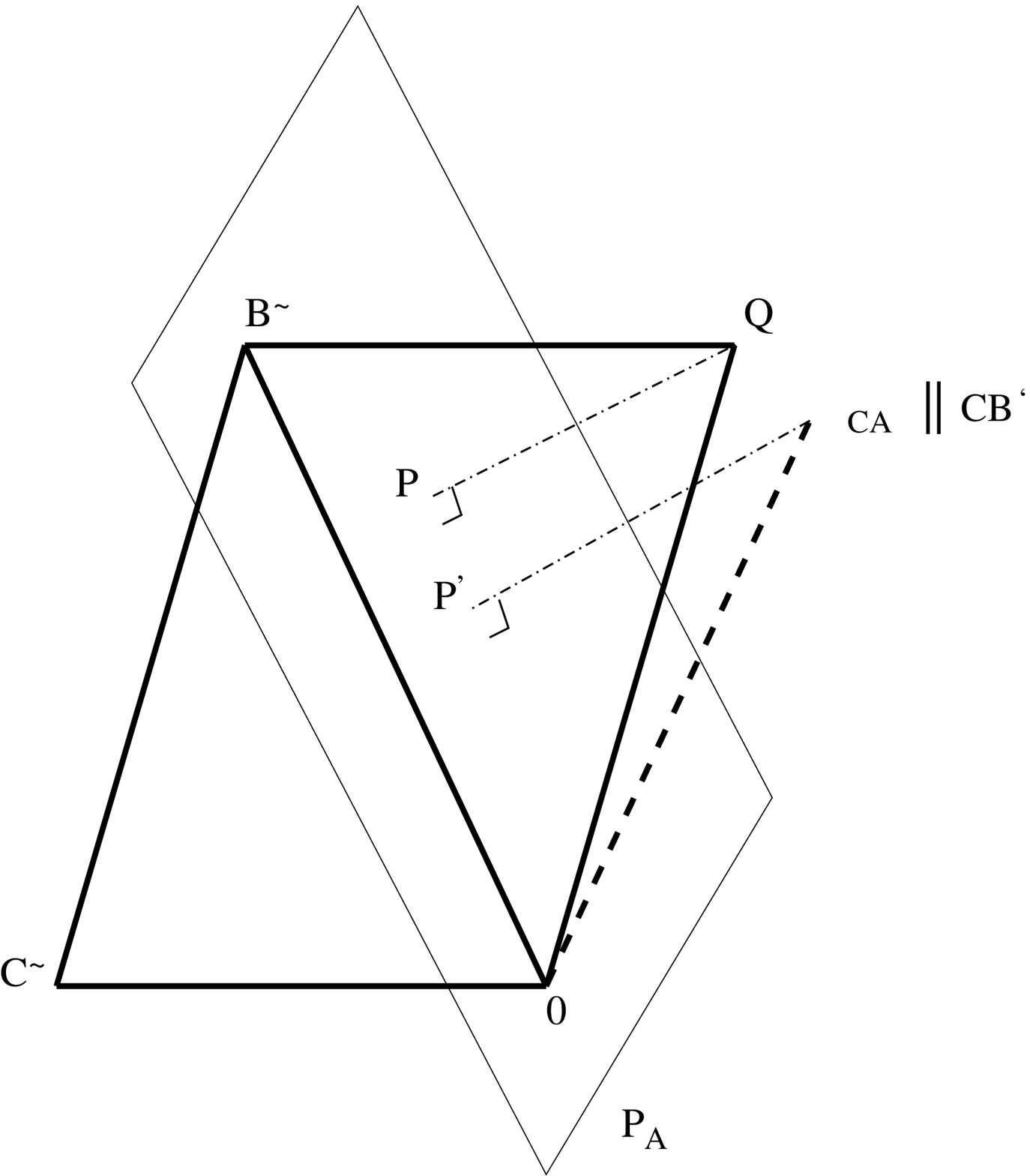,height=7cm}}

\noindent \textit{Claim 3: }$$d(C,B) \ge \frac{20}{21}\frac{1}{
({\cal C}+6) \sqrt{1 + {\tan}^2 \delta}}.$$\\
The triangle $\tb P \tbc$ is a right \tr at $P$, so
using claim 1 and claim 2 we have
\begin{equation*} \begin{split}
|\tbc - \tb |^2 =& |\tb -P|^2 + |\tbc -P|^2\\
\le&  \left(1 +  {\tan}^2 \delta\right)\cdot |\tbc
-P|^2 \\
\le& \left(1 +  {\tan}^2 \delta\right)\cdot (6+{\cal
C})^2\cdot d(C,A)^2\cdot |\tbc -0|^2. \end{split}
 \end{equation*}
The vector $\tbc -\tb$ is just $0- \tc$, the length of
which is $d(A,C)$, and the vector $\tbc -0$ is $ \tb
-\tc$. So $$d(A,C) \le \sqrt {(1 +  {\tan}^2 \delta)}
(6+{\cal C})\cdot d(C,A)|\tb -\tc|, $$ and
$$\frac{1}{(6+ {\cal C}) \sqrt{1 + {\tan}^2
\delta}}\le |\tb - \tc|.$$ Using standard estimates - see
Corollary \ref{CorA4} - we obtain $ |\tb - \tc| \le \frac{11}{10}
d(C,B)$, and the
proposition follows. \Kasten\\[0.2cm]

\section{Application of Proposition  \ref{Prop31} to $\vertg$.}\label{s5}

Fix $g$ in $G$. Let $C$ and $A$ be points on Weinstein's average $N$ with $d(C,A)<
0.15$ joint by a minimizing \geo $\gamma$ in $M$. Suppose that
$\exp_C(v)=\exp_A(w)=:B$ for vertical vectors $v \in \vertg_C$ and
$w \in \vertg_A$ of lengths less than 0.5. In this section we will
apply Proposition \ref{Prop31} to the geodesic triangle given by
the above three points of $M$ and $P_A=\vertg_A,
P_C=\vertg_C$. We will do so in Proposition \ref{Prop45}\\[0.2cm]
To this aim, first we will estimate the constants $\delta$ and
${\cal C}$ of Proposition  \ref{Prop31} in this specific case. As
always our estimates will hold for $\e<\frac{1}{20000}$.\\Roughly
speaking, the constant $\delta$ - which measures how much the
angle between $\dot{CA}=\dot\gamma(0)$ and $\vertg_C$ deviates
from $\frac{\pi}{2}$ - will be determined by using the fact that
$N$ is $C^1$-close to $N_g$, so that the shortest geodesic
$\gamma$ between $C$ and $A$ is ``nearly tangent'' to the
distribution $Hor^g$.\\ Bounding the constant ${\cal C}$ - which
measures how the angle between $\vertg_A$ and $\vertg_C$ depends
on $d(A,C)$ - will be easier, noticing that both spaces are
parallel translations of normal spaces to $N_g$, which is a
submanifold with bounded
second fundamental form.\\[0.2cm]
Since we have
$$\measuredangle\left(\dgamma(0), \vertg_C \right)=
\frac{\pi}{2} - \measuredangle\left(\dgamma(0), \horg_C \right),$$
to determine $\delta$ we just have to estimate the angle
$$\alpha:= \measuredangle\left(\dgamma(0), \horg_C \right).$$ We
already introduced the geodesic $\gamma(t)$ from $C$ to $A$, which
we assume to be parametrized to arc length. We now consider the
curve $\pi(t):=\png \circ \gamma(t)$ in $N_g$. We can lift the
curve $\pi$ to a curve $\vg \circ \pi$ in $N$ connecting $C$ and
$A$; we will call $c(t)$ the  parametrization to arc length of
this lift.\\ $\np$ will denote the connection induced on $\nu N_g$
by the Levi-Civita connection $\nabla$ of $M$, and $
^{\perp}_{\pii} \parap$ applied to some $\xi \in \nu_{\pi(t)}N_g$
will denote its $\np$-parallel transport from $\pi(t)$ to $\pi(0)$
along $\pi$. (The superscript ``$b$'' stands for ``backwards'' and
is a reminder that we are parallel translating to the initial
point of the curve
$\pi$.)\\
Further we will need $$r:=100\e + \frac{L(\gamma)}{2}
\ge \sup_{t}\{d(\gamma(t),N_g)\} \text{    and   }
f(r):=\cos(r) -\frac{3}{2}\sin(r).$$
Notice that $r<0.08$ due to our restrictions on $\e$
and $d(C,A)$.\\
\centerline{\epsfig{file=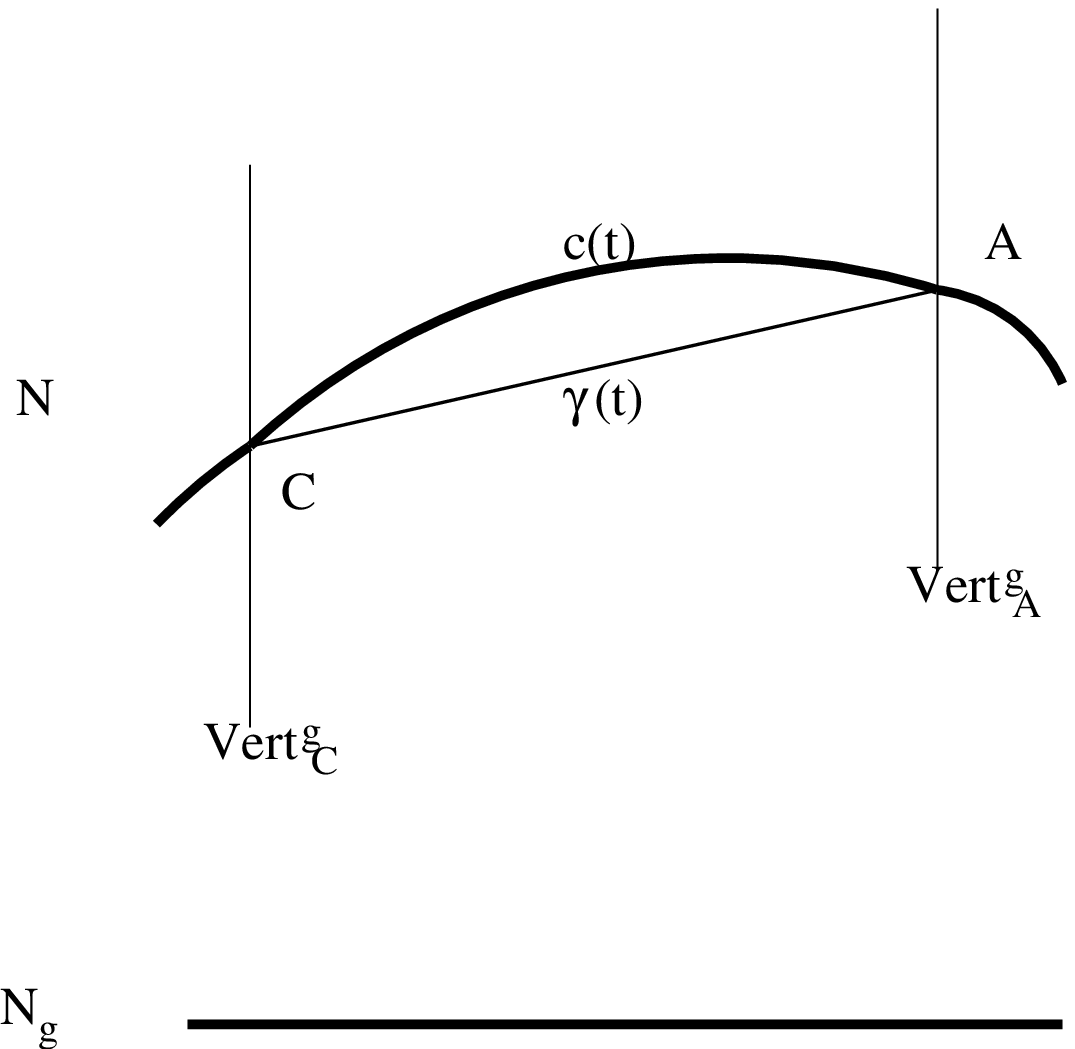,height=7cm}}

\noindent Using the fact that $c$ is a curve in $N$ and $N$ is
$C^1$-close to $N_g$, in Appendix B we will show that the section
 $\tilde{c}:=\expi_{N_g}(c(t))$ of $\nu N_g$ along
$\pi$ is ``approximately parallel''. This will allow
us to bound from above the ``distance'' between its
endpoints as follows:
\begin{Prop} \label{Prop41} $$\left|\exp_{N_g}^{-1} C-
_{\pii}^{\perp} \parap (\exp_{N_g}^{-1} A) \right| \le L(\gamma)
\frac{3150 \e}{f(r)}.$$ \end{Prop} Using the fact that $\gamma$ is
a geodesic and our bound on the extrinsic curvature of $N_g$, in
Appendix C we will show that the section
$\tgamma:=\expi_{N_g}(\gamma(t))$ of $\nu N_g$ along $\pi$
approximately ``grows at a constant rate''. Since its covariant
derivative at zero depends on $\alpha$, we will be able to
estimate the ``distance'' between its endpoints  - which are also
the endpoints of $\tilde{c}$ - in terms in $\alpha$. We will
obtain:
\begin{Prop}\label{Prop42}$$\left|\expi_{N_g}(C) -  _{\pi^b}^{\perp}
\parap \expi_{N_g}(A)\right| \ge L(\gamma)\left[
\frac{99}{100}\sin \left(\alpha -\frac{\epsilon}{4} \right) - 500
\e - 3r - \frac{8}{3}L(\gamma)\left(r+\frac{r + \frac{3}{2}}{f(r)}
\right)\right].$$ \end{Prop} Comparison of Proposition
\ref{Prop41} and  Proposition \ref{Prop42} gives
$$ \frac{3150 \e}{f(r)} \ge  \frac{99}{100}\sin
\left(\alpha -\frac{\epsilon}{4} \right) - 500 \e - 3r
- \frac{8}{3}L(\gamma)\left(r+\frac{r +
\frac{3}{2}}{f(r)}\right).$$
Recall that $r= 100\e+ \frac{L(\gamma)}{2}$. If
$L(\gamma)$ and $\epsilon$ are small enough one can
solve the above inequality for $\alpha$. With our
restriction $\e <\frac{1}{20000}$ this can be done whenever
$L(\gamma)<0.1$. One obtains
$$\delta(\epsilon,L(\gamma)):= \frac {\epsilon}{4} +
\arcsin \left\{ \frac{100}{99}\left[\frac{3150 \e}{f(r)} + 500 \e
+ 3r + \frac{8}{3}L(\gamma)\left(r+\frac{r +
\frac{3}{2}}{f(r)}\right)\right]\right\}> \alpha.$$ We can now
state the main results of this section. First we determine the
constant $\delta$ of Proposition  \ref{Prop31} in our setting.
\begin{Prop}\label{Prop43} Let $C,A$ be points in $N$ and $\gamma$
the shortest \geo in $M$ from $C$ to $A$. Assume
$\e<\frac{1}{20000}$ and $L(\gamma)< 0.1$. Then
$\delta(\epsilon,L(\gamma))$ is well defined and
$$\measuredangle\left(\horg_C,
\dot{\gamma}(0)\right)=\alpha< \delta
\left(\epsilon,L(\gamma) \right).$$
Therefore $$\measuredangle\left(\vertg_C,
\dgamma(0)\right) \ge \frac{\pi}{2}- \delta
\left(\epsilon,L(\gamma) \right)$$ and for symmetry
reasons
 $$\measuredangle\left(\vertg_A, -\dgamma(L(\gamma))
\right) \ge \frac{\pi}{2}- \delta\left(\epsilon,L(\gamma)
\right).$$ \end{Prop} To determine the constant $\cal C$ we only
need Lemma \ref{LemC3}:
\begin{Prop}\label{Prop44} Let $C,A$ and $\gamma$ be as above and
assume $L(\gamma)<0.1$. Then
$$d(\vertg_C, _{\gammai}\parap \vertg_A) \le
2L(\gamma).$$ \end{Prop} \pf By Lemma \ref{LemC3} we have
$d(\vertg_C,\, _{\gammai}\parap \vertg_A) \le  \arcsin
[L(\gamma)(r+ \frac{r+\frac{3}{2}}{f(r)})]$, where $r= 100\e+
\frac{L(\gamma)}{2}$. For the above values of $\e$ and
${L(\gamma)}$ this last expression is bounded
above by $2L(\gamma)$. \Kasten\\[0.2cm]
Now making use of the estimates in the last two propositions we
can apply Proposition  \ref{Prop31}.
\begin{Prop}\label{Prop45} Fix $g \in G$.
Let $C,A$ be points in $N$ \st $d(A,C)<0.1$  and suppose that
$\exp_C(v)=\exp_A(w)=:B$ for vertical vectors $v\in \vertg_C, w\in
\vertg_A$. Then
$$|v|,|w|\ge \frac{1}{9}\frac{1}{\sqrt{1+\tan^2
\left(\delta(\e, d(A,C)\right)}}.$$ \end{Prop}
\pf If $|v|\ge 0.5$
than the estimate for $|v|$
clearly holds, as the right hand side is $\le \frac{1}{9}$. So we assume $|v|=d(B,C)<0.5$. \\
Since $d(B,N_g)<0.5+100\e<1$ and $(M,N_g)$ is a gentle pair, the
triangle $ABC$ lies in an open subset of $M$ with the
properties\begin{quote}
i) the sectional curvature lies between -1 and 1\\
ii) the injectivity radius at each point is at least 1.
\end{quote} Therefore we are in the situation of Proposition
\ref{Prop31}. Setting $P_C=\vertg_C$ and $P_A=\vertg_A$ in the
statement of Proposition  \ref{Prop31},  Proposition \ref{Prop43}
and
Proposition  \ref{Prop44} allow us to choose:\\
 $$\delta=\delta \left( \epsilon, d(A,C)\right) \text{
 and  } {\cal C}=2.$$
Therefore, since $\frac{10}{11}(\frac{1}{2+6})>\frac{1}{9}$, we
obtain
$$|v| \ge  \frac{1}{9}\frac{1}{\sqrt{1+\tan^2
\left(\delta(\e, d(A,C)\right)}}.$$
The statement for $|w|$ follows exactly in the same
way. \Kasten

\section{Estimates on tubular neighborhoods of $N_g$
on which $\vg$ is injective}\label{s6}

In this section we will finally apply the results of Section
\ref{s4} and Section \ref{s5}, which were summarized in
Proposition \ref{Prop45}, to show that $\exp_{N_g}(\nu
N_g)_{0.05}$ is a tubular \nbhd of $N$ on which $\vg$ is
injective. We will also bound from below the size of $\cap_{g \in
G} \exp_N (Vert^g)_{0.05}$ (where the 2-form $\int_g\omega_g$ is
defined).
\begin{Prop}\label{Prop51} If $\e <\frac{1}{20000}$ the map $$\vg:
\exp_{N_g}(\nu N_g)_{0.05} \rightarrow \exp_N (Vert^g)_{0.05}$$ is
a diffeomorphism. \end{Prop} \pf From the definition of $\vg$ it
is clear that it is enough to show the injectivity of $$\exp_N:
(\vertg)_{0.05} \rightarrow \exp_N(\vertg)_{0.05}.$$ Let $A,C\in
N$ and $v\in \vertg_C, w\in \vertg_A$ be vectors of length
strictly less than $0.05$. We suppose that  $\exp_C(v)=\exp_A(w)$
and argue by
contradiction.\\
Clearly $d(A,C)< 0.1$. We can apply Proposition  \ref{Prop45},
which implies $|v|,|w|\ge
\frac{1}{9}\frac{1}{\sqrt{1+\tan^2(\delta(\epsilon, d(A,C))}}.$
Since the function $\delta(\epsilon, L)$ increases with $L$ we
have $$|v|,|w|\ge \frac{1}{9}\frac{1}{\sqrt{1+\tan^2
\left(\delta(\epsilon, 0.1 )\right)}}.$$ For $\e <\frac{1}{20000}$
the above function is larger than $0.05$, so we have a
contradiction.\\ So $\exp_C(v)\neq \exp_A(w)$ and the above map is
injective. \Kasten\\[0.2cm]
For each $L \le 0.05$ we want to estimate the radius of a tubular
\nbhd of $N$ contained in $\cap_{g \in G} \exp_N (Vert^g)_{L}$.
This will be used in Section \ref{s7} to determine where $\int
\omega_g$ is non-degenerate, so that one can apply Moser's trick
there. As a by-product, the proposition below will also give us an
estimate of the size of the neighborhood in which $\int_g
\omega_g$ is defined.
 \begin{Prop}\label{Prop52} For $L \le 0.05$ and $\e <\frac{1}{20000} $, using
the notation $$\re_L:= \sin (L)\cos \left(\delta(\epsilon, 2L)+
2L^2\right)$$ we have
$$\exp_N(\nu N)_{\re_L}\subset  \cap_{g \in G} \exp_N (Vert^g)_{L}.$$
 \end{Prop}
\rem The function $\re_{0.05}$ decreases with $\e$ and
assumes the value $0.039...$ at $\e=0$ and the value
$0.027...$ when $\e=\frac{1}{20000}$.\\[0.4cm]
To prove the proposition we will consider again  \geo
triangles:
\begin{Lem}\label{Lem51} Let $ABC$ be a \geo \tr lying in $\tnng$
\st $d(A,B)\le d(C,B) =: L < 0.05$. Let $\gamma$ denote the angle
at $C$, and suppose that $\gamma \in [\frac{\pi}{2}- \td,
\frac{\pi}{2}+ \td]$. Then $$d(A,B) \ge \sin(L) \cos \left(\td
+2L^2 \right).$$ \end{Lem} \pf Denote by $\alpha, \beta$ the
angles at $A$ and $B$ respectively, and denote further by
$\alpha', \beta',\gamma'$ the angles of the Alexandrov triangle in
$S^2$ corresponding to $ABC$ (i.e. the triangle in $S^2$ having
the same side lengths as $ABC$). By [Kl, remark 2.7.5] we have
$\sin(d(A,B))=\sin(d(C,B))\frac{\sin(\gamma')}{\sin(\alpha')}\ge
\sin(d(C,B))\sin(\gamma')$.\\By Toponogov's theorem (see [Kl])
$\gamma' \ge \gamma$. On the other hand, using the bound $L^2$ for
the area of the Alexandrov triangles in $S^2$ and $H^2$
corresponding to $ABC$, we have $\gamma'-\gamma \le 2L^2$ (see
proof of Lemma \ref{LemA1}). So $\gamma' \in [\frac{\pi}{2} -
\tilde{\delta},
 \frac{\pi}{2} + \tilde{\delta}+ 2L^2]$. Altogether
this gives
$$d(A,B) \ge \sin\left( d(A,B) \right) \ge
\sin(d(C,B))\sin(\gamma')\ge \sin(L)\sin \left(
 \frac{\pi}{2} + \tilde{\delta}+ 2L^2 \right).$$
\Kasten\\[0.2cm]
\noindent Now we want to apply the Lemma \ref{Lem51} to our case
of interest:
\begin{Lem}\label{Lem52} Let $C\in N$ and $B=\exp_C(w)$ for some
$w\in \vertg_C$ of length $L<0.05$, and assume as usual
$\e<\frac{1}{20000}$. Then $$d(B,N)  \ge
 \sin(L) \cos \left(\delta(\epsilon,2L)
+2L^2\right)=R^{\epsilon}_L.$$ Here the function $\delta$ is as in
Section \ref{s5}.
\end{Lem} \pf Let $A$ be the closest point in $N$ to $B$. Clearly
$d(A,B) \le d(C,B)=L$, so the shortest geodesic $\gamma$ from $C$
to $A$ has length $L(\gamma) \le 2L$. By Proposition \ref{Prop43}
we have
 $$\measuredangle \left( \dgamma(0),\vertg_C \right)
\ge \frac{\pi}{2}- \delta \left(\epsilon,L(\gamma)
\right)
 \ge \frac{\pi}{2}- \delta(\epsilon, 2L).$$
So, since $w\in \vertg_C$, $$\measuredangle( \dgamma(0),w) \in
\left[\frac{\pi}{2}- \delta(\epsilon,2L)\;,\; \frac{\pi}{2}+
\delta(\epsilon,2L) \right].$$ Using the fact that, for any $g\in
G$, the triangle $ABC$ lies in $\exp_{N_g}(\nu N_g)_1$, the lemma
follows using Lemma \ref{Lem51} with $\tilde{\delta}=\delta(\e,
2L)$.\Kasten\\[0.2cm]
\textit {Proof of Proposition  \ref{Prop52}:} For any $g\in G$ and
positive number $L<0.05$, by Lemma \ref{Lem52} each point $B \in
\partial(\exp_N(\vertg)_L)$ has distance at least $\sin(L) \cos
(\delta(\e, 2L) +2L^2)=R^{\e}_L$ from $N$. Therefore $tub(\re_L)$
lies in $\exp_N(\vertg)_L$, and since this holds for all $g$ we
are done. \Kasten

\section{Conclusion of the proof of the Main Theorem}\label{s7}

 In sections \ref{s3}-\ref{s6},
making use of the Riemannian structure of $M$, we showed that the
two-form $\int_g \omega_g$ is well-defined in the neighborhood
$\cap_{g \in G} \exp_N (Vert^g)_{0.05}$ of $N$ (Recall that
$\omega_g:=(\vg^{-1})^* \omega$ was defined in the
Introduction). In this section we will focus on the symplectic
structure of $M$ and conclude the proof of the Main Theorem, as
outlined in Part II of Section \ref{s13}.\\[0.4cm]
First we will show that $\iog$ is a symplectic form on a suitably
defined neighborhood $\te$ of $N$. Then it will easily follow that
the convex linear combination $\ot:=\omega +t(\iog -\omega)$
consists of symplectic forms.\\ As we saw in the introduction,
$[\omega]=[\iog]\in H^2(\te,\RR)$, so we can apply Moser's trick. The
main step consists of constructing canonically a primitive
$\alpha$  of small maximum norm for the two-form $\frac{d}{dt}
\ot$.\\ Comparing the size of the resulting Moser \vf with the
size of $\te$ we will determine an $\e$ for which the existence of
an
isotropic average of the $N_g$'s is ensured.\\[0.4cm]
In this section we require $L<0.05$. Notice that the estimates of
Section \ref{s3} hold for such $L$. We start by requiring
$\e < \frac{1}{20000}$ and introduce the abbreviation
$$\de_L:=4L+4100\e$$ for the upper bound obtained  in
Proposition  \ref{Prop24} on  $\exp_N(\vertg)_{L}$.
\\

\subsection{Symplectic forms in $\te$}\label{s71}

In Section \ref{s3} we estimated the difference between $(\vg)_*X$
and $
\para X$. This lemma does the same for $\vgi$.
\begin{Lem}\label{Lem61} Let $q\in \partial \tnnl$
and $X \in T_qM$ a unit vector. Then $$\left|(\vgi)_*X -\para X
\right|\le \frac{\de_L}{1-\de_L}.$$ Furthermore,
$$\frac{1}{1+\de_L} \le \left|(\vgi)_*X \right|\le
\frac{1}{1-\de_L}.$$ \end{Lem} \pf Let $p:=\vgi(q)$. By
Proposition  \ref{Prop24}, for any vector $Z\in T_pM$ we have

$$\frac{|\vg_*(Z)|}{1+ \de_L} \le \left|Z \right| \le
\frac{|\vg_*(Z)|}{1- \de_L}.$$
The second statement of this lemma follows setting
$Z=(\vgi)_*X$.\\
Choosing instead $Z=(\vgi)_*X - \para X \in T_pM$ and applying
once more Proposition  \ref{Prop24} gives
$$\left| (\vgi)_*X - \para X \right| \le \frac{|X -
(\vg)_* \para X|}{1 -\de_L} \le
\frac{\de_L}{1-\de_L}.$$ \Kasten\\[0.2cm]
\noindent Since $(\vgi)_*X$ is close to $\para X$ and since our
assumption on $\nabla \omega$ allows us to control to which extent
$\omega$ is invariant under parallel translation we are able to
show that $\omega$ and $\omega_g=(\vgi)^*\omega$ are close to each
other:
\begin{Lem}\label{Lem62} Let $X,Y$ be unit tangent vectors at  $q \in \exp_N(\vertg)_{L}$. Then
 $$
|\left({\og -\omega})(X,Y)\right| \le \frac{\de_L}{1-\de_L}\left(
\frac{\de_L}{1 -\de_L}+2\right)+2L+100\e.$$
\end{Lem}
\pf Denoting $p:=\vgi(q)$ we have
 \begin{equation*} \begin{split}
(\og -\omega)_q(X,Y)=& \omega_p \Big( (\vgi)_*X \;,\;
(\vgi)_*Y \Big) - \omega_q(X,Y)\\
=&  \omega_p \Big(\para X - \left[ (\vgi)_*X  -
\para X \right] \;,\;\para Y - \left[ (\vgi)_*Y  -
\para Y \right] \Big) - \omega_q(X,Y)\\
=& \omega_p \Big( (\vgi)_*X  - \para X \;,\; (\vgi)_*Y
 - \para Y \Big)\\
&+ \omega_p \Big(\;\para X  \;,\;  (\vgi)_*Y  - \para Y \Big) +
\omega_p \Big(  (\vgi)_*X  - \para X  \;,\; \para Y \Big)\\
&+ \omega_p(\para X, \para Y)-\omega_q(X,Y).
\end{split} \end{equation*}
Now since ``$\para$'' is the parallel translation along a curve of
length $<2L+100\e$ (see Section \ref{s3}) and $|\nabla \omega|<1$
we have $ \omega_p(\para X, \para Y)-\omega_q(X,Y)<2L+100\e$ and
using Lemma \ref{Lem61}
we are done. \Kasten\\[0.2cm]
\noindent Since the symplectic form $\omega$ is compatible with
the metric and the $\omega_g$'s are close to $\omega$ we obtain
the non-degeneracy of $\omega_t$ for $L$ and $\epsilon$ small
enough.

\begin{Cor}\label{Cor61}
Let $X$ be a unit tangent vector at  $q \in \cap_{g\in
G}\exp_N(\vertg)_{L}$. Then for all $t\in [0,1]$
$$\omega_t(X,IX) \ge1- \left[ \frac{\de_L}{1- \de_L}\left( \frac{\de_L}{1- \de_L} +2
\right)+2L+100\e \right].$$
\end{Cor}
 \pf By definition
\begin{equation*}\begin{split}
\omega_t(X,IX) = \omega(X,IX)+t\cdot \int_g
(\omega_g-\omega)(X,IX).
\end{split}\end{equation*}
The first term is equal to 1 because $\omega$ is almost-K\"ahler,
the norm of the second one is estimated using Lemma \ref{Lem62}.
 \Kasten\\[0.2cm]

\noindent \rem The right hand side of Corollary \ref{Cor61} is surely
positive if $\de_L\le 0.1$. We set \footnote{This choice of $\lep$ will
 allow us to obtain good numerical estimates in
Section \ref{s74}.}
$$\lep:= \frac{0.1 - 4100 \e}{4} $$ and require $\e < \frac{1}{70000}$.
We obtain immediately:
 \begin{Prop}\label{Prop61} On
$$
\te:=
\exp_N(\nu N)_{R^{\epsilon}_{L^{\epsilon}}} \subset \cap_{g\in G} \exp_N(\vertg)_{\lep}$$
the convex linear combination  $\ot:=\omega +t(\iog -\omega)$
consists of symplectic forms. \end{Prop}
\rem
Recall that the function $R^{\epsilon}_L$ was defined in
Proposition  \ref{Prop52}. See Section \ref{s74} for a graph of
$R^{\epsilon}_{L^{\epsilon}}$, a function of $\epsilon$.
\\

\subsection{The construction of the primitive of
$\frac{d}{dt} \ot$}\label{s72}

We want to construct canonically a primitive $\alpha$ of $\frac{d}{dt} \ot =
\int \omega_g -\omega$ on  $\cap_{g\in G} \exp_N(\vertg)_{0.05}$.
We first recall the following fact, which is a slight
modification of [Ca, Chapter III].\\
Let $N$ be a submanifold of a Riemannian \mfd $M$, and let $E
\rightarrow N$ be a subbundle of $TM|_N \rightarrow N$ \st
$E\oplus TN=TM|_N$. Furthermore let $\tilde{U}$ be a fiber-wise
convex \nbhd of the
 zero
section of $E \rightarrow N$ \st $\exp: \tilde{U}
\rightarrow U \subset M$ be a diffeomorphism. Denote
by $\pi: U \rightarrow N$ the projection along the
slices given by exponentiating the fibers of $E$ and
by $i: N \hookrightarrow M$ the inclusion.\\ Then
there is an operator $Q: \Omega^{\bullet}(U)
\rightarrow  \Omega^{\bullet -1}(U)$ \st
$$\text{Id}-(i \circ \pi)^*=dQ+Qd:\;
\Omega^{\bullet}(U) \rightarrow \Omega^{\bullet}(U).$$ A concrete
example is given by considering $\rho_t:U \rightarrow U,\;
\exp_q(v) \mapsto \exp_q(tv)$ and ${w_t}|_{\rho_t(p)}:=
\frac{d}{ds}|_{s=t}\rho_s(p)$. Then $$Qf:=\int_0^1Q_tf \text{dt},
\;\;\;\; Q_tf:= \rho_t^*(i_{w_t}f)$$ gives an operator with the
above
property.\\
Note that for a 2-form $\omega$ evaluated on $X\in
T_pM$ we have
\begin{equation*}\begin{split} \left|(Q_t \omega)_p
X\right| =&\left| \omega_p \left(w_t|_{\rho_t(p)},
{\rho_t}_*(X)\right)\right| \\ \le& |\tilde{\omega}_p|_{\text{op}}
\cdot d(p, \pi(p)) \cdot | {\rho_t}_*(X)|\;\;\;
\;\;\;\;\;\;\;\;\;\;\;\;\;\;
\;\;\;\;\;\;\;\;\;\;\;\;\;\;\;\;\;\;\;\;\; (\bigstar)
\end{split} \end{equation*}
where $ |\tilde{\omega}_p|_{\text{op}}$ is the
operator norm of $\tilde{\omega}_p:T_pM \rightarrow
T_p^*M$ and the inner product on $T_p^*M$ is induced
by the one on $T_pM$.\\[0.4cm]
\noindent For each $g$ in $G$ we want to construct a canonical
primitive $\alpha$ of $\og - \omega$ on $\exp_N(\vertg)_{0.05}$
 We do so in
two steps:\\[0.4cm]
\noindent \begin{large} \textbf{Step I  } \end{large}  We apply
the above procedure to the \vb $\vertg \rightarrow N$ to obtain an
operator $\qgn$ \st $$\text{Id}- (\pgn)^* \circ (i_N)^*=d \qgn
+\qgn d$$ for all differential forms on $\exp_N(\vertg)_{0.05}$.
Since $N$ is isotropic with respect to $\og$ and $\og - \omega$ is
closed we have
$$\og -\omega =d \qgn(\og -\omega)
+(\pgn)^*(i_N)^*(-\omega).$$

\noindent \begin{large} \textbf{Step
II } \end{large}  Now we apply the procedure to the
\vb $\nu N_g \rightarrow N_g$ to get an operator
$\qng$ on differential forms on $\exp_{N_g}(\nu
N_g)_{100\e}$. Since $N_g$ is isotropic with respect
to $\omega$ we have
$$ \omega = d \qng \omega,$$ so we found a primitive
of $\omega$ on $\exp_{N_g}(\nu N_g)_{100\e}$.\\
Since $N\subset \exp_{N_g}(\nu N_g)_{100\e}$ the
1-form $ \bg := i_N^*(\qng \omega)$ on $N$ is a
well defined primitive of $i_N^*\omega$.\\[0.4cm]
Summing up these two steps we see that
$$ \ag:= \qgn(\og - \omega) - (\pgn)^* \bg$$ is a
primitive of $\og - \omega$ on
$\exp_N(\vertg)_{0.05}$.\\[0.2cm]
So clearly $\alpha:=\int_g \ag$ is a primitive of
$\frac{d}{dt}\omega_t= \iog - \omega$ on  $\cap_{g\in G} \exp_N(\vertg)_{0.05}$.

\subsection{Estimates on the primitive of
$\frac{d}{dt} \ot$}\label{s73}

In this section we will estimate the $C^0$-norm of the one-form $\alpha$ constructed
in Section \ref{s72}.\\[0.4cm]
\noindent \begin{large} \textbf{Step II } \end{large} We will first estimate
the norm of $ \bg := i_N^*(\qng \omega)$ using $(\bigstar)$ and
then the norm of $(\pgn)^* \bg$.
\begin{Lem}\label{Lem63} If $p \in \exp_{N_g}(\nu N_g)_{100 \e}$
and $X\in T_pM$ is a unit vector, then for any $t\in [0,1]$
$$\left|\rngt_*X\right| \le \frac{5}{4}.$$ \end{Lem}
\pf Let $L:=d(p,N_g)<100 \e < \frac{1}{700}$ and write
$X=J(L)+K(L)$, where $J$ and $K$ are $N_g$ Jacobi
fields along the unit speed geodesic $\gamma_p$ from
$p':=\png(p)$ to $p$ \st $J(0)$ vanishes, $J'(0)$ is
normal to $N_g$, and $K$ is a strong $N_g$ Jacobi
field (see the remark in Section \ref{s32}).\\
$J(t)$ is the variational \vf of a variation
$f_s(t)=\exp_{p'}(tv(s))$ where the $v(s)$'s are unit normal
vectors at $p'$. Therefore $${\rngt}_* J(L)=\frac{d}{ds} \Big|_0
\left[ \rngt \circ \exp_{p'}(Lv(s)) \right] = \frac{d}{ds} \Big|_0
\left[\exp_{p'}(tLv(s)) \right]=J(tL).$$ Using Lemma \ref{Lem24}
we have on one hand $|LJ'(0)| \le (1+\sigma(L))|J(L)|$ and on the
other hand $|J(tL)|(1-\sigma(tL)) \le tL|J'(0)|$ where $
\sigma(x):=\frac{\sinh(x)-x}{\sin(x)}.$ So
$$|J(tL)| \le
t \frac{1+\sigma(L)}{1-\sigma(tL)}|J(L)| \le
\frac{21}{20}t|J(L)|.$$ Similarly we have ${\rngt}_* K(L)=K(tL)$.
Using Lemma \ref{Lem26} we have $|K(0)| \le (1+
\frac{9}{5}L)|K(L)|$ and $|K(tL)| \le \frac{|K(0)|}{1-
\frac{9}{5}tL}$, therefore $|K(tL)| \le \frac{1+
\frac{9}{5}L}{1- \frac{9}{5}tL}|K(L)| \le \frac{21}{20} |K(L)|.$\\
 Altogether we
have
\begin{equation*} \begin{split} \left|{\rngt}_*X\right|^2 =&
\left| J(tL)+K(tL)
\right|^2\\
\le&  \left| J(tL) \right|^2 + \left| K(tL) \right|^2
+ 2\cdot \frac{9}{5} tL\left| J(tL) \right|\left|
K(tL) \right| \\
\le&  (\frac{21}{20})^2 \left[  \left| J(L) \right|^2 + \left| K(L)
\right|^2 +  \frac{18}{5}L  \left| J(L) \right|\left|
K(L) \right| \right]\\
\le& ( \frac{21}{20})^2 \left[  \left| J(L) + K(L) \right|^2 +
\frac{36}{5}L  \left| J(L) \right|\left| K(L) \right|
\right]\\
\le& ( \frac{21}{20})^2 \left[ 1+\frac{36}{5}\cdot 1.1^2L \right]
\le \frac{5}{4}
\end{split} \end{equation*}
where in the first and third inequality we used Lemma \ref{Lem29}
and in the fourth in addition we used Lemma
\ref{Lem210}.\Kasten\\[0.2cm]
\begin{Cor}\label{Cor62} The one-form $\beta^g$ on $N$ satisfies $$
|\bg| < 125 \e.$$  \end{Cor} \pf At any point $p\in N \subset
\exp_{N_g}(\nu N_g)_{100 \e}$, using $(\bigstar)$ ,
$|\tilde{\omega}|_{\text{op}}=1$ and Lemma \ref{Lem63}, we have
$|(\qng \omega)_p|< 125 \e $. Clearly $|(\qng \omega)_p| \ge
|i_N^*(\qng \omega)_p|=|\bg_p|.$
\Kasten\\[0.2cm]
Now we would like to estimate $(\pgn)_*X$ for a unit
tangent vector $X$. Since $\pgn={(\rho_N^g)}_0$ we
prove a stronger statement that will be used again
later. Recall that we assume $L\le 0.05$.
\begin{Lem}\label{Lem64} If $q\in \exp_N(\vertg)_{L} $ and $X\in
T_qM$ is a unit vector then for any $t \in [0,1]$ we have
$$|\rgnt_*X| \le 1.5 \frac{1 +\de_L}{1-\de_L}.$$ \end{Lem}
\pf Using Lemma \ref{Lem61} we have
$$\left|{\rgnt}_*X \right| \le (1+ \de_L) \left| (\vgi)_*{\rgnt}_*X
\right|.$$Clearly $\rgnt \circ \vg = \vg \circ \rngt$, since - up
to exponentiating - $\vg$ maps $\nu N_g$ to $\vertg$, and $\rngt$
and $\rgnt$ are just rescaling of the
respective fibers by a factor of $t$.\\
If we reproduce the proof of  Lemma \ref{Lem63} requiring $p$ to
lie in $\exp_{N_g} (\nu N_g)_{L}$ we obtain
\footnote{Since $L<0.05$ now we have to replace the constant $\frac{21}{20}$ in that proof by the constant $\frac{6}{5}$.}
$|{\rngt}_*Y|<1.5$ for unit vectors $Y$ at $p$. Using this and
Lemma \ref{Lem61} respectively we have
$$ \left|{\rngt}_*(\vgi)_*X \right| \le 1.5 \left|(\vgi)_*X
\right| \;\;\;\text{  and  }\;\;\; \left|(\vgi)_*X \right| \le
\frac{1}{1-\de_L}.$$
Altogether this proves the lemma. \Kasten\\[0.2cm]

\begin{Cor}\label{Cor63} On $ \exp_N(\vertg)_{L} $  we have
$$|(\pgn)^* \bg | \le  200 \e  \frac{1 +\de_L}{1-\de_L}.$$
\end{Cor}
\pf This is clear from $|((\pgn)^*
\beta^g)X|=|\beta^g((\pgn)_*X)|$, Corollary \ref{Cor62} and
Lemma \ref{Lem64}.\Kasten\\[0.2cm]
\begin{large} \textbf{Step I  } \end{large}  Now we
estimate $|\qgn(\og -\omega)|$. This is easily achieved using
Lemma \ref{Lem62} and Lemma \ref{Lem64} to estimate the quantities
involved in $(\bigstar)$:
\begin{Cor} \label{Cor64}  For $q\in \partial \tnnl$
we have $$\left|\qgn(\og -\omega)_q\right| \le 1.5 \frac{1+
\de_L}{1- \de_L}\cdot L \cdot \left[ \frac{\de_L}{1-\de_L}\left(
\frac{\de_L}{1 -\de_L}+2\right) +2L+100\e\right].$$ \end{Cor} \rem By Proposition
5.2
$d(q,N) \ge \re_L$. Furthermore, when $\e<\frac{1}{70000}$ and $L<
0.05$, one can show that $\re_L \ge \frac{2}{3}L$. So $L\le \frac{3}{2}\,
d(q,N).$\\
[0.4cm] Now finally using Corollary \ref{Cor63} and Corollary
\ref{Cor64} we can estimate the norm of $\alpha:=\int_g \ag$:
\begin{Prop}\label{Prop62}  Assuming $L<0.05$ at $q\in \cap_g\exp_N(\vertg)_{L} $
we have
\begin{eqnarray*}
\left|\alpha_q \right| &\le&
 1.5 \frac{1+ \de_L}{1- \de_L}\cdot \frac{3}{2} d(q,N) \cdot
\left[\frac{\de_L}{1-\de_L}\left( \frac{\de_L}{1 -\de_L}+2\right)+2L+100\e\right]\\
&+& 200 \e  \frac{1 +\de_L}{1-\de_L}.
\end{eqnarray*}
\end{Prop}


\subsection{The end of the proof of the Main Theorem}\label{s74}
Proposition  \ref{Prop61} showed that the Moser vector field
$v_t:=-\tilde{\omega_t}^{-1} \alpha$ is well-defined on $\te
\subset \cap_{g\in G} \exp_N(\vertg)_{\lep}$.
 Recalling that
$D^{\epsilon}_{L^{\epsilon}}=0.1$, Corollary \ref{Cor61}
immediately implies
\begin{Cor} \label{Cor65}  At $q\in \te \subset \cap_g\exp_N(\vertg)_{L^{\e}} $
we have
$$\left|(\widetilde{\ot})_q^{-1} \right| _{\text{op}}
\le \frac{1}{1- \left[ \frac{\de_{L^{\e}}}{1-\de_{L^{\e}}}( \frac{\de_{L^{\e}}}{1
-\de_{L^{\e}}}+2) +2L^{\e} +100\e \right]} \le 1.53.$$
\end{Cor}

\noindent From Corollary \ref{Cor65} and Proposition \ref{Prop62}
we obtain:
\begin{Prop}\label{Prop63} For all $t\in [0,1]$ and $q \in \te$
$$\left|(v_t)_q \right| \le
\left|(\widetilde{\ot})_q^{-1}
\right|_{\text{op}}\cdot \left| \alpha_q \right|
\le
1.45\,d(q,N) + 374 \e.$$
\end{Prop}

\noindent Let $\gamma(t)$ be an integral curve of the
time-dependent vector field $v_t$ on $\te$ \st $p:=\gamma(0)\in
N$. Where $d(\,\cdot\,,p)$ is differentiable, its gradient has
unit length. So $
\frac{d}{dt}d(\gamma(t),p) \le |\dgamma(t)|$.\\
By Proposition  \ref{Prop63} we have $ |\dgamma(t)| \le 1.45\,
d(\gamma(t),p)+ 374 \e$. So altogether
$$ \frac{d}{dt}d \left(\gamma(t),p \right) \le 1.45\, d
\left(\gamma(t),p \right)+ 374 \e.$$ The solution of the ODE
$\dot{s}(t)=As(t)+B$ satisfying $s(0)=0$ is
$\frac{B}{A}(e^{At}-1)$. Hence, if the integral curve $\gamma$ is
well defined at time 1, we have
$$ d(\gamma(1),N)\le  d(\gamma(1),p) \le \frac{374\e}{1.45}\left( e^{1.45}-1 \right) \le 842 \e.$$
Let us denote by $\rho_1$ the time-1 flow of the time dependent
vector field $v_t$, so that $\rho_1^{-1}$ is the time-1 flow of
$-v_{1-t}$.
Since by definition $\te:=
\exp_N(\nu N)_{R^{\epsilon}_{L^{\epsilon}}}$ the submanifold $L:=\rho_1^{-1}(N)$
will surely be well defined if
$$842 \e < R^{\epsilon}_{L^{\epsilon}}.$$
This is always the case since $\e < \frac{1}{70000}$.\\

\includegraphics{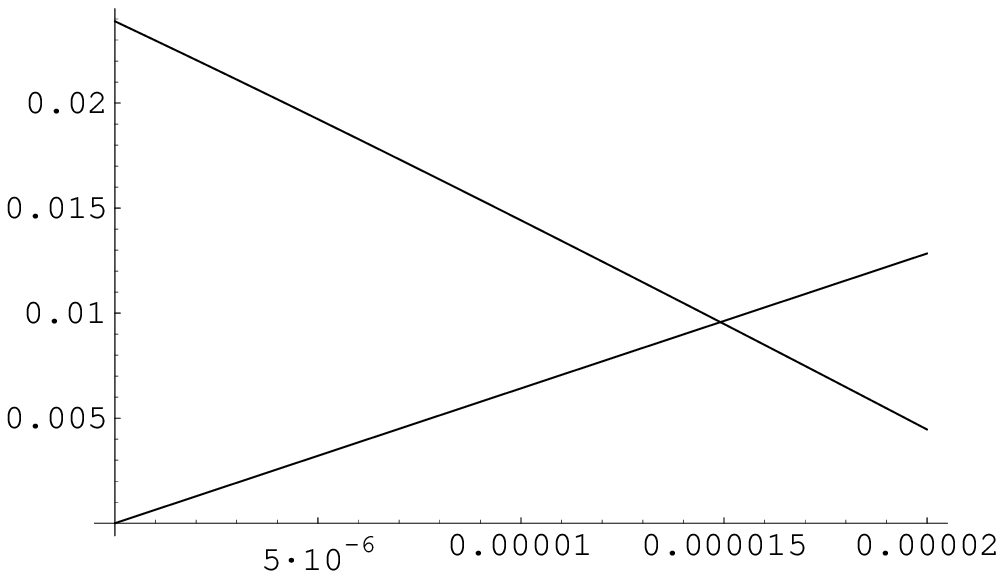}

\noindent \textsc{Graphs of $842\e$ (increasing) and
$R^{\epsilon}_{L^{\epsilon}}$
(decreasing).} \\[0.4cm]
The estimate for $d_0(N_g,L)$ is obtained using
$d_0(N,L)<842 \e$ and $d_0(N_g,N)<100 \e$.
The proof of the Main Theorem is now complete.\\

\noindent \textbf{Remark: } In the Main Theorem we assumed that  $|\nabla
\omega|<1$. Let us now consider the case that  $|\nabla \omega|$
is bigger than one. Then the statement of the  Main Theorem still
holds verbatim if one makes the bound on  $\epsilon$ smaller, as
follows. The bound on $|\nabla \omega|$ enters our proof directly
only in Lemma \ref{Lem62}; if $|\nabla \omega|$ is bigger than
one, the right hand side of that lemma should read $|\left({\og
-\omega})(X,Y)\right| \le \frac{\de_L}{1-\de_L}\left(
\frac{\de_L}{1 -\de_L}+2\right)+|\nabla \omega|(2L+100\e)$
instead. Similarly, the quantity ``$2L+100 \e$'' appearing in
Corollary \ref{Cor61}, Corollary \ref{Cor64} and Proposition
\ref{Prop62} should be multiplied by $|\nabla \omega|$. Now let us
assume that $\epsilon<\frac{1}{|\nabla \omega|} \frac{1}{70000}$
and let us replace $L^{\e}$ everywhere by
$\tilde{L}^{\e}:=\frac{\frac{0.1}{|\nabla \omega|}-4100 \e}{4}.$
Then the bounds on $\left|(\widetilde{\ot})_q^{-1} \right|
_{\text{op}}$ and $|(v_t)_q|$ given in Corollary \ref{Cor65} and
Proposition \ref{Prop63} still hold, and our isotropic average $L$
will be well defined if $842 \e <
R^{\epsilon}_{\tilde{L}^{\epsilon}}$. This is satisfied for
$\epsilon$ small enough, since $R^{\epsilon}_{\tilde{L}^{\epsilon}}$ is a continuous function and $R^{0}_{\tilde{L}^{0}}$ is positive.

\section{Remarks on the Main Theorem}\label{s8}

\textbf{Remark 1: } \textit{Is the isotropic average
$L$ $C^1$ close to the $N_g$'s?}\\
The main shortcoming in our Main Theorem is surely
the lack of an estimate on the $C^1$ distance
$d_1(N_g,L)$.\\ To bound $d_1(N_g,L)$
it is enough to estimate the distance between tangent
spaces $T_pL$ and $T_{\rho_1(p)}N$. Indeed, this would
allow us to estimate the distance between $T_pL$ and
$T_{\png(p)}N_g$, using which - when $\e$ is small
enough - one can conclude that $\png: L \rightarrow
N_g$ is a diffeomorphism and give the desired bound on
the $C^1$ distance.\\
Using local coordinates and standard theorems about
ODEs it is possible to estimate the distance between
tangent spaces $T_pL$ and $T_{\rho_1(p)}N$ provided
one has a bound on the covariant derivative of the
Moser vector field, for which one would have to
estimate $\nabla(\tilde{\omega_t})^{-1}$. To do so one
should be able to bound expressions like
$\nabla_Y((\vgi)_*X)$ for parallel vector fields $X$
along some curve.\\This does not seem to be possible
without more information on the extrinsic geometry of
$N$. We recall that it is not known whether the
average $N$ forms a gentle pair with $M$, see Remark
6.1 in [We]. We are currently trying to improve Weinstein's theorem so that one obtains a gentle average.
\\
[0.4cm]
\textbf{Remark 2: } \textit{The case of isotropic
$N$.}\\
Unfortunately, if the Weinstein average $N$ happens to be already
isotropic with respect to $\omega$, our construction will
generally provide an isotropic average $L$ different from $N$.
Indeed, while Step I of Section \ref{s72} always gives a one-form
vanishing on points of $N$, Step II does not, even if $N$ is
isotropic for $\omega$.\\
The procedure outlined in Remark 3, on the other hand,
would produce $N$ as the isotropic average, but in
that case the upper bound for $\e$ would depend on the
geometry of $N$.\\[0,4cm]
\textbf{Remark 3: } \textit {Averaging of symplectic
and coisotropic submanifolds.}\\
The averaging of $C^1$-close gentle symplectic submanifolds of an almost-K\"ahler
manifold is a much simpler task than for isotropic
submanifolds. The reason is that $C^1$-small perturbations of
symplectic manifolds are symplectic again and one  can simply
apply
Weinstein's averaging procedure ([We, Thm. 2.3]).\\
Unfortunately our construction does not allow to
average coisotropic submanifolds. In our proof we were able to canonically construct a
primitive of $\iog -\omega$ using the fact that the
$N_g$'s are isotropic with respect to $\omega$. In
the case that the $N_g$ are not isotropic
  it is still possible to construct
canonically a primitive, following Step I of our
construction and making use of the primitive
$d^*(\triangle^{-1} i^*_N(\og -\omega))$ of $i^*_N(\og
-\omega)$ (but the upper bound on its norm would
depend on the geometry of $N$).\\
Nevertheless, our construction fails in the
coisotropic case, since the fact that $N$ is
coisotropic for all $\og$'s does not imply that it is
for their average $\iog$.\\[0.4cm]

 \section{An application to Hamiltonian actions}\label{s9}
As a simple application of our Main Theorem
we apply Theorem 2 to almost invariant isotropic submanifolds
of a Hamiltonian G-space and deduce some information
about their images under the moment map.\\

\noindent We start by recalling some basic definitions (see [Ca]):
consider an action of a Lie group $G$ on a symplectic manifold
$(M,\omega)$ by symplectomorphisms. A moment map for the action is
a map $J:M \rightarrow \cg^*$ such that for all $v \in \cg$ we
have $ \omega(v_M,\cdot)=d \langle J,v \rangle$ and which is
equivariant with respect to the $G$ action on $M$ and the
coadjoint action of $G$ on
 $ \cg$. Here $v_M$ is the vector field on $M$ given
by $v$ via the
infinitesimal action. An action admitting a moment map is called Hamiltonian action.\\

\noindent This simple lemma is a counterpart to ([Ch], Prop 1.3).
\begin{Lem}\label{Lemc1}
Let the compact connected Lie group $G$ act on the symplectic
manifold $(M,\omega)$ with moment map $J$. Let $L$ be a connected
isotropic submanifold of $(M,\omega)$ which is invariant under the
group action. Then $L\subset J^{-1}(\mu)$ where $\mu$ is a fixed
point of the coadjoint action.
\end{Lem}
\textit{Proof: } Let $X \in T_x L$.  For each $v\in \cg$ we have
$$d_x\langle J,v \rangle X= \omega(v_M(x),X)=0,$$
since both $v_M(x)$ and $X$ are tangent to the isotropic
submanifold $L$. Therefore every component of the moment map is
constant along $L$, so $L\subset J^{-1}(\mu)$ for
some $\mu \in \cg^*$.\\
Now let $x_0\in L$ and $G\cdot x_0 \subset L$ the orbit through
$x_0$. Then from the equivariance of $J$ it follows that for all
$g$ we have $\mu=J(g \cdot x_0)=g\cdot J(x_0)=g\cdot \mu$, so
$\mu$ is a fixed point of the coadjoint action. \Kasten

\noindent Now we apply the lemma above to the case where $L$ is
almost invariant.

\begin{Cor}\label{Corc1} Let the compact Lie group $G$ act on the symplectic
manifold $(M,\omega)$ with moment map $J:M \rightarrow \cg^*$.
Suppose $M$ is endowed with a $G$-invariant compatible Riemannian
metric so that the Levi-Civita connection satisfies $|\nabla
\omega|<1$. If a
connected isotropic submanifold
$L\subset M$ satisfies:\\
i) $(M,L)$ is a gentle pair\\
ii) $d_1(L, g\cdot L)<\epsilon <\frac{1}{70000}$ for all $g \in G$\\
then $J(L)$ lies the ball of radius $1000\epsilon \cdot C$ about
a fixed point $\mu$ of the coadjoint action.
Here $\cg^*$ is endowed with any inner product and
$C:=Max\{|v_M|:
v\in \cg \text{ has unit
length}\}$.\\

\end{Cor}

\noindent \textit{Proof: } By Theorem 2 there exists an isotropic
submanifold $L'$ invariant under the $G$ action with
$d_0(L,L')<1000\epsilon $. By Lemma \ref{Lemc1} $L$ lies the some
fiber $J^{-1}(\mu)$ of the moment map, where $\mu$ is a fixed
point of the coadjoint action. We will show that
$|J(p)-\mu|<1000\epsilon \cdot  C$ for all $p\in L$. \\
Let $p'$ a closest point to $p$ in $L'$. The shortest geodesic   $\gamma$ from $p$ to $p'$, which
we choose to be
defined on the interval $[0,1]$, has length $<1000\e$.
Therefore for any unit-length
$v\in \cg$ (with respect to the inner product induced on $\cg$ by
its dual) we have
\begin{eqnarray*}\langle J(p)-\mu,v \rangle &=&
 \int_0^1 \langle dJ(\gamma(t))\dot{\gamma}(t),v
\rangle dt\\
 &=& \int_0^1 d\langle J,v\rangle \dot{\gamma}(t) dt\\
  &=& \int_0^1 \omega(v_M,\dot\gamma(t)) dt.
\end{eqnarray*}
Since for all $t$ we have $|\omega(v_M,\dot\gamma(t))|\le |v_ M |
\cdot|\dot\gamma(t)|\le 1000\epsilon \cdot C $ we are done.
\Kasten

\appendix
\section{The estimates of Proposition  \ref{Prop31}}

Here we will prove the estimates used in the proof of the
Proposition  \ref{Prop31}. See Section \ref{s4} for the
notation.\\[0.4cm]
Applying a series of lemmas we will show first that
$$\measuredangle(P-0,P'-0)\le 2d(C,A).$$ Then we will
show $$\measuredangle( (\expi_A)_* \dot{CB}\,,\, \tb -\tc) <
3d(A,C),$$ from which very easily it will follow that
$$\measuredangle(_{CA}\parap \dot{CB}\,,\,\tbc-0)<4 d(C,A).$$ This
will conclude the proof of claim 2 in the proof of Proposition
 \ref{Prop31}.\\ Finally using standard arguments we will obtain the
estimate used in claim 3, namely
$$ |\tb - \tc| \le \frac{11}{10}
d(C,B).$$\\[0.4cm]
We begin by stating an easy result without
proof:\\[0.2cm]
\textbf{Fact } \textit{ Let $c:[0,1] \rightarrow \RR^n$ is a
differentiable curve \st  $|\dot{c}(s) -\dot{c}(0)|\le
\varepsilon$ for all $s \in [0,1]$. Then
$$\big|[c(1)-c(0)]-\dot{c}(0) \big| \le \varepsilon.$$
If in addition ${\cal
C} \le |\dot{c}(s)| \le {\cal D}$ for all $s \in [0,1]$,
 then} $$\cos \Big(\measuredangle[c(1) -c(0)\,,\,
\dot{c}(0)] \Big) \ge \frac{{\cal C}^2
-\frac{\varepsilon^2}{2}}{{\cal D}^2}.$$
Now we introduce the following notation: $\sigma(s)$
will be the shortest \geo from $C$ to $B$, defined on
$[0,1]$, and $\tsigma(s)$ will be its lift to $T_AM$
via $\exp_A$. For any $s\in [0,1]$, $f_s(t)=\exp_A(t
\tsigma(s))$ will be the shortest \geo from $A$ to
$\sigma(s)$. Notice that, since $d(A,C)+d(C,B) \le
0.7$, we have $|\tsigma(s)|\le 0.7$.\\[0.4cm]
For a fixed $s$, consider the \geo \tr $A \sigma(s) C$
in $M$, and denote by $\alpha,\beta_s$ and $\gamma$
the angles at $A, \sigma(s)$ and $C$ respectively.
\begin{Lem}\label{LemA1} $$\text{If
}\frac{8}{7}\frac{d(C,A)}{\sin(70^{\circ})} \le | \tsigma(s)|
\text{, then } \sin(\beta_s) \le \frac{8}{7}\frac{d(C,A)}{ |
\tsigma(s)|} \le \sin(70^{\circ}).$$ \end{Lem}
 \rem The lemma says that if $\sigma(s)$ is far enough
from $A$, the direction of the geodesic $\sigma$ there
does not deviate too much from being  radial \wrt $A$.
This will allow us to estimate the radial and
orthogonal-to-radial components of
$\dot{\sigma}(s)$ in Lemma \ref{LemA3} and Lemma \ref{LemA5}.\\[0.2cm]
\pf Consider the Alexandrov triangle in $S^2$ corresponding to
$A\sigma(s)C$, i.e. a triangle in $S^2$ with the same side lengths
as $A\sigma(s)C$ and denote its angles by $\alpha', \beta_s'$ and
$\gamma'$. Consider also the corresponding triangle in standard
hyperbolic space $H^2$ and denote its angles by  $\alpha'',
\beta_s''$ and $\gamma''$. The sum of the angles of the triangle
in $S^2$ deviates from $180^{\circ}$ by the area of the triangle,
which is bounded by $0.7 \cdot 0.15 < 10^{\circ}$ (see [BK] 6.7
and 6.7.1.), so that the sum is at most $190^{\circ}$. For the
same reason the sum of the angles of the triangle in $H^2$ is at
least
 $170^{\circ}$. Since for the sectional curvature
$\kappa$ we have $|\kappa|\le 1$, by Toponogov's
Theorem (see Thm. 2.7.6 and Thm. 2.7.12 in [Kl]) we
have $\beta_s'' \le \beta_s \le \beta'_s $ (and
similarly for $\alpha, \gamma$) and by the above
$\beta_s' - \beta_s''<20^{\circ}$.\\[0.2cm]
\textit{Case 1: }$\beta_s\le 90^{\circ}$.\\ By the law
of sines (see [Kl, Remark 2.7.5]) we have
$$\sin \beta_s'=\sin(d(A,C))\cdot \frac{\sin
\gamma'}{\sin\, d(\sigma(s),A)} \le
\frac{d(A,C)}{\sin\,d(\sigma(s),A)} \le \frac{8}{7}
\frac{d(A,C)}{\,d(\sigma(s),A)},$$ where we used
$\frac{sin x}{x} > \frac{7}{8}$ for $x \in [0,0.8]$.
By our assumption this quantity is bounded above by
$\sin(70^{\circ})$, so $\beta_s'\in [0,70^{\circ}]$ or
  $\beta_s'\in [110^{\circ},180^{\circ}].$ Since
$\beta_s' - \beta_s<20^{\circ}$ and $\beta_s\le
90^{\circ}$, we must have  $\beta_s'\in
[0,70^{\circ}]$. Therefore $\beta_s \le \beta_s'$
implies $\sin \beta_s \le \sin \beta_s'$, which we
already bounded above.\\[0.2cm]
\textit{Case 2: }$\beta_s > 90^{\circ}$.\\
One has to proceed analogously, but comparing with
triangles in $H^2$. \Kasten\\[0.2cm]
We now state a general fact about the exponential map:
\begin{Lem}\label{LemA2} If $\gamma$ is a geodesic parametrized to
arc-length and $W \in T_{\gamma(0)}M$, then for $t<0.7$ $$ \left|
_{\gamma}\parap ( d_{t \dot{\gamma}(0)} \exp_{\gamma(0)}) W - W
\right| \le \frac{\sinh(t)-t}{t}|W|\le \frac{t^2}{5}$$ and
$$ \left(1-\frac{t^2}{5}\right)|W| \le |\exp_*W|\le
\left(1+\frac{t^2}{5}\right)|W|.$$ \end{Lem}
 \pf The unique \jf
along $\gamma$ \st $ J(0)=0, J'(0)=W$ is given by $J(t)=(d_{t
\dot{\gamma}(0)} \exp_{\gamma(0)})(tW)$ (see [Jo, Cor. 4.2.2]).
The estimate $\frac{\sinh(t)-t}{t}$ follows from [BK 6.3.8iii].
This expression is bounded above by $\frac{t^2}{5}$ when $t<0.7$.
The second estimate follows trivially from the first one. We
prefer to use these estimates rather than more standard ones (see
 [BK 6.4.1])
in order to keep the form of later estimates  more concise.\Kasten\\[0.2cm]
Using Lemma \ref{LemA1} we obtain a refinement of Lemma
\ref{LemA2} for $\dot{\sigma}(s)=(d_{\tsigma(s)} \exp_A)
\dtsigma(s)$:
\begin{Lem}\label{LemA3} For all $s \in [0,1]$ $$\big| _{f_s}
\parap \dsigma(s) - \dtsigma(s) \big| \le
\frac{1}{5}d(C,A)|\dot{\sigma}(s)|.$$ \end{Lem} \pf By the Gauss
Lemma, the decomposition into radial (denoted by $^R$) and \otr
(denoted by $^{\perp}$) components of $\dot{\sigma}(s)=
(d_{\tsigma(s)} \exp_A) \dtsigma(s)\in T_{\sigma(s)}M$ and
$\dtsigma(s) \in T_{\tsigma(s)}(T_AM)$ is preserved by
$d_{\tsigma(s)} \exp_A$. Hence $_{f_s} \parap \dsigma(s)^R=
_{f_s}\parap  (d_{\tsigma(s)} \exp_A) [\dtsigma(s)^R]=
\dtsigma(s)^R \in T_AM$. Introducing the notation $D:=
|\tsigma(s)| < 0.7$ this implies
$$\big| _{f_s}\parap \dsigma(s) - \dtsigma(s)\big|
=\big|_{f_s} \parap \dsigma(s)^{\perp} - \dtsigma(s)^{\perp}\big|
\le \frac{D^2}{5} |\dtsigma(s)^{\perp}|$$ by Lemma \ref{LemA2}.
From this we obtain $ |\dtsigma(s)^{\perp}| \le
\frac{|\dsigma(s)^{\perp}|}{ 1-  \frac{D^2}{5} } \le
\frac{5}{4}|\dsigma^{\perp}(s)|$. So $| _{f_s}\parap \dsigma(s) -
\dtsigma(s)| \le
\frac{D^2}{4}|\dot{\sigma}^{\perp}(s)|$.\\[0.2cm]
\textit{Case 1: } $\frac{8}{7}
\frac{d(C,A)}{\sin(70^{\circ})} \le D$.\\
Using the notation of Lemma \ref{LemA1} we have
$|\dsigma(s)^R|=|\cos \beta_s|\cdot |\sigma(s)|$, so
$|\dot{\sigma}(s)^{\perp}|= \sin ( \beta_s)|\dsigma(s)| $.
Applying the same lemma we obtain $|\dot{\sigma}(s)^{\perp}|\le
\frac{8}{7} \frac{d(C,A)}{D} |\dot{\sigma}(s)|$. Therefore
$\frac{D^2}{4} |\dsigma(s)^{\perp}| \le \frac{D^2}{4}\frac{8}{7}
\frac{d(C,A)}{D} |\dot{\sigma}(s)| \le \frac{1}{5}
 d(C,A)|\dot{\sigma}(s)|$.\\[0.2cm]
\textit{Case 2: } $\frac{8}{7}
\frac{d(C,A)}{\sin(70^{\circ})} > D$.\\
We simply use $ |\dsigma(s)^{\perp}|\le  |\dsigma(s)|$
to estimate $$\frac{D^2}{5} |\dsigma(s)^{\perp}|\le
\left( \frac{8}{7} \frac{d(C,A)}{\sin(70^{\circ})}
\right)^2 \frac{1}{5} |\dsigma(s)| \le \frac{3}{10}
d(C,A)^2 |\dsigma(s)|,$$ which is less than  the
estimate of \textit{Case 1} since $d(C,A)\le
0.15$.\Kasten\\[0.2cm]
We want to apply the above Fact to the curve
$\tsigma:[0,1] \rightarrow T_AM$.
We begin by  determining the constant $\varepsilon$.
\begin{Lem}\label{LemA4} For all $s \in [0,1]$ $$ | \dtsigma(s) -
\dtsigma(0)|\le d(C,A)\cdot d(C,B). $$ \end{Lem} \pf We will show
that $\dtsigma(s)\approx \,\, _{f_s}
\parap \dsigma(s) \approx \,\,  _{f_0} \parap
\dsigma(0) \approx \,\dtsigma(0)$, where in the first
and last relation tangent spaces are identified by
parallel translation along $f_s$ and $f_0$
respectively. \\ The estimates needed for the first
and last relation are given by Lemma \ref{LemA3} since $|\dot{\sigma}(s)|=d(C,B)$ for all $s$.\\
For the second relation we use the fact that $\sigma$ is a
geodesic and [BK, 6.2.1]. We see that $| _{f_s} \parap \dsigma(s)
- _{f_0}
\parap \dsigma(0)|$ is bounded above by the product of
$|\dsigma(0)|=d(C,B)$ and of the area spanned by the triangle
$\sigma(s)AC$ - which is surely less than $0.6 \,d(C,A)$ by [BK,
6.7.1] - and we are
done.\Kasten\\[0.2cm]
Now we determine the constants ${\cal C}$ and ${\cal D}$. Recall
that ${|\dsigma(s)|}=d(C,B)$, so that the constants we obtain are
really independent of $s$.
\begin{Lem}\label{LemA5} For all $s \in [0,1]$
$$\frac{|\dsigma(s)|}{\sqrt{1 + d(C,A)^2}} \le
|\dtsigma(s)| \le\frac{|\dsigma(s)|}{\sqrt{1 - d(C,A)^2}} .$$
\end{Lem} \rem Since we need estimates involving $d(C,A)$, the
classical Jacobi-field estimate  $  |\dtsigma(s)| \le
\frac{D}{\sin(D)}|\dsigma(s)|$ or $  |\dtsigma(s)| \le
(1-\frac{D^2}{5})^{-1}|\dsigma(s)|$ (given by Lemma \ref{LemA2})
are not good enough for us. Here $D:=d(\sigma(s),A)$. Furthermore,
making use of Lemma \ref{LemA1}, in Lemma \ref{LemA3} we bounded
above $| _{f_s}
\parap \dsigma(s) - \dtsigma(s) |$, which at once implies
  $ |\dtsigma(s)| \le
(1+\frac{d(C,A)}{5}) |\dsigma(s)|$. This however is also not
sufficient, because it would only allow us to bound the left hand
side of Corollary \ref{CorA2} by a multiple of $\sqrt{d(C,A)}$.
Instead, we will use Lemma \ref{LemA1} to bound directly
$\frac{|\dot{\sigma}(s)|}{|\dot{\tilde{\sigma}}(s)|}$ in order to obtain a
bound of the form $[1+o(d(C,A)^2)]|\dot{\sigma}(s)|$ in Lemma \ref{LemA5}.
\\[0.2cm]
\pf We fix $s$ and adopt the notation
$D:=|\tsigma(s)|=d(\sigma(s),A)<0.7$.\\
\textit{Case 1: } $D< \frac{8}{7}
\frac{d(C,A)}{\sin(70^{\circ})} $.\\
In Case 2 of the proof of Lemma \ref{LemA3} we showed that  $|
_{f_s} \parap  \dsigma(s) -   \dtsigma(s)| \le \frac{3}{10}
d(C,A)^2 |\dsigma(s)|.$ So
$$|\dsigma(s)|\left(1 -  \frac{3}{10} d(C,A)^2 \right)
\le |\dtsigma(s)| \le |\dsigma(s)|\left(1 +
\frac{3}{10} d(C,A)^2 \right).$$\\[0.2cm]
\textit{Case 2: } $D \ge \frac{8}{7}
\frac{d(C,A)}{\sin(70^{\circ})}$.\\
In view of the remark above we will use simple Jacobi field
estimates only for the \otr components of $\dsigma(s)$ and
$\dtsigma(s)$, which we can bound above using Lemma \ref{LemA1},
whereas for the radial components we just have to notice that they
have the same length.\\Recall that $f_s(t)=\exp_A(t \tsigma(s))$
is a variation of geodesics emanating from $A$. We denote by
$\tilde{f}_s(t)=t \tsigma(s)$ its lift to $T_AM$. Fixing $s\in
[0,1]$, we denote by $J(t)$ and $\tj(t)$ the Jacobi-fields along
$f_s$ and $\tilde{f}_s$ arising from the variations  $f_s(t)$ and
$\tilde{f}_s(t)$. So $J(t)=(d_{t \tsigma(s)}\exp_A)t\dtsigma(s)$
and $\tj(t)=t \dtsigma(s)$ have the same initial covariant
derivative $E:= \dtsigma(s)$.\\ We decompose $E\in T_AM$ into its
components parallel and orthogonal to $\tilde{\sigma}(s)$ as
$E=E^R+E^{\perp}$. $J^R(t):=(d_{t \tsigma(s)}\exp_A)tE^R$ and
$J^{\perp}(t):=(d_{t \tsigma(s)}\exp_A)tE^{\perp}$ are Jacobi
fields, since $\exp_A$ maps lines through $0 \in T_AM$ to
geodesics through $A \in M$. They both vanish at zero and their
initial covariant derivatives are $E^R$ and $E^{\perp}$
respectively, so by the uniqueness of Jacobi fields to given
initial data we have $J=J^R+J^{\perp}$. The Gauss Lemma implies
that this is the decomposition of $J$ into radial and
orthogonal-to-radial components.\\[0.2cm]
To show the second inequality in the statement of the
lemma we want to bound below
$$ \frac{|\dsigma(s)|^2}{|\dtsigma(s)|^2}=
\frac{|J(1)|^2}{|\tj(1)|^2}= \frac{|E^R|^2 +
|J^{\perp}(1)|^2}{|E|^2}.$$ Applying Lemma \ref{LemA2} to
$J(1)^{\perp}=(d_{\sigma(s)}\exp_A)\dtsigma(s)^{\perp}$ we obtain
$|J(1)^{\perp}| \ge (1- \frac{D^2}{5})|E^{\perp}|$, so
$$ \frac{|\dsigma(s)|^2}{|\dtsigma(s)|^2} \ge
\frac{|E^R|^2 + (1- \frac{D^2}{5})^2|E^{\perp}|^2}{|E|^2}=
1+\left(- \frac{2}{5}\frac{D^2}{3}+\frac{D^4}{25}
\right)\frac{|E^{\perp}|^2} {|E|^2} .\;\;(\blacktriangledown)$$
Now we bound from above  $ \frac{|E^{\perp}|^2} {|E|^2}$. We
already saw that $|E^{\perp}| \le
\frac{|J^{\perp}(1)|}{1-\frac{D^2}{5}}$, and by Lemma \ref{LemA2}
$|\dsigma(s)| \le (1+ \frac{D^2}{5}) |\dtsigma(s)|$. Both
$\frac{1}{1-\frac{D^2}{5}}$ and $1+ \frac{D^2}{5}$ are bounded by
$\frac{8}{7}$ since $D< 0.7$. Now Lemma \ref{LemA1} allows us to
relate $|J^{\perp}(1)|$ and $|\dsigma(s)| \in T_{\sigma(s)}M$.
Namely, we have  $|J^{\perp}(1)| = \sin (\beta_s)|\dsigma(s)| \le
\frac{8}{7} \frac{d(C,A)}{D} |\dsigma(s)|.$ The last three
estimates give $|E^{\perp}| \le (\frac{8}{7})^3
\frac{d(C,A)}{D}|E|.  $ Substituting into $(\blacktriangledown)$
gives
$$ \frac{|\dsigma(s)|^2}{|\dtsigma(s)|^2} \ge  1 +
\left(\frac{8}{7} \right)^6  \left(-
\frac{2}{5}\frac{D^2}{3}+\frac{D^4}{25} \right) \frac{
d(C,A)^2}{D^2} \ge
  1 -   \left(\frac{8}{7} \right)^6 \frac{2}{5}
d(C,A)^2.$$
To take care of the  first inequality in the statement
of the lemma we show that
$$ \frac{|\dsigma(s)|^2}{|\dtsigma(s)|^2} \le  1 +
\left(\frac{8}{7} \right)^6 \frac{3}{7} d(C,A)^2$$ by repeating
the above proof and using the estimate
$|J^{\perp}(1)| \le (1+ \frac{D^2}{5})|E^{\perp}|$.\\[0.2cm]
To finish the proof we have to compare the estimates obtained in
\textit{Case 1} and \textit{Case 2}. To do so notice that
$(\frac{8}{7} )^6 \frac{3}{7} <1$, that $(1 +  \frac{3}{10}
d(C,A)^2 ) \le \frac{1}{\sqrt{1 - d(A,C)^2}}$ and  $
\frac{1}{\sqrt{1 +  d(A,C)^2}}\le (1 -  \frac{3}{10} d(C,A)^2 )$.
\Kasten\\[0.2cm]
Now finally we can apply the above Fact to the curve
$\tilde{\sigma}:[0,1]\rightarrow T_AM$. The first statement of the
Fact allows us to prove
\begin{Cor}\label{CorA1} $$ \measuredangle(P-0,P'-0) \le 2d(C,A).$$
\end{Cor}
\pf We first want to bound $|P'-P|$ from above and
$|P'-0|$ from below.\\Since $P'$ and $P$ are the
closest points in $P_A$ to $_{CA}\parap \dot{CB}$
and $Q$ respectively,
\begin{equation*} \begin{split} \left|P-P' \right|
\le& \left|(Q-0)- _{CA}\parap \dot{CB}\right| \\ \le&
\left|(Q-0) - (\expi_A)_* \dot{CB} \right| + \left|
(\expi_A)_* \dot{CB} -  _{CA}\parap \dot{CB}\right| \\
\le& d(C,A) |\dot{CB}|+\frac{d(C,A)^2}{5}|\dot{CB}|.
\end{split}\end{equation*}
In the last inequality we used
$Q-0=\tb-\tc=\tilde{\sigma}(1)-\tilde{\sigma}(0)$ and $(\expi_A)_*
\dot{CB}=\dtsigma(0)$ to apply the first statement of the above
Fact (with $\varepsilon$ given by Lemma \ref{LemA4}) for the first
term and Lemma \ref{LemA2} for the
second term.\\
On the other hand we have
\begin{equation*} \begin{split} |P'-0 | =& |\dot{CB}
|\cdot \cos \big(\measuredangle(_{CA}\parap \dot{CB},
P_A) \big)  \\ \ge& |\dot{CB}| \cdot \cos(\theta) \\
\ge& |\dot{CB}|\sqrt{1 - \theta^2} \\ \ge&
|\dot{CB}|\sqrt{1 - {\cal C}^2 d(C,A)^2}.
\end{split}\end{equation*}
Therefore we have $$\sin \big(\measuredangle(P'-0,P-0) \big) \le
\frac{ |P'-P|}{ |P'-0|}\le \frac{\frac{d(C,A)}{5}+1}{\sqrt{1 -
{\cal C}^2 d(C,A)^2}} d(C,A).$$ So, using the restrictions ${\cal
C} \le 2$ and $d(C,A)<0.15$, and using $\frac{\sin x}{x} \ge
\frac{7}{8}$ for $x\in [0,0.8]$, we obtain $$
\measuredangle(P'-0,P'-P) \le \frac{8}{7} \sin
\bigl(\measuredangle(P'-0,P'-P) \bigr)  \le
2d(C,A).$$\Kasten\\[0.2cm]
The second statement of the above Fact delivers
\begin{Cor}\label{CorA2}  $$ \measuredangle((\expi_A)_* \dot{CB},
\tb -\tc) < 3 d(C,A).$$ \end{Cor} \pf Applying again the above
Fact to the curve $\tsigma: [0,1] \rightarrow T_AM$ we want to
obtain an estimate for $\alpha :=\measuredangle ( \tsigma(1)-
\tsigma(0)\; , \; \dtsigma(0)) = \measuredangle( \tb -\tc \;,\;
(\expi_A)_* \dot{CB}).$ Lemma \ref{LemA4} and Lemma \ref{LemA5}
deliver the estimates
$$ \varepsilon = d(C,A)\cdot d(C,B), \;\;\;\;  {\cal
C}=  \frac{d(C,B)}{\sqrt{1 + d(C,A)^2}}  \;\;\;\text{ and } {\cal
D}= \frac{d(C,B)}{\sqrt{1 - d(C,A)^2}} .$$ Therefore, using the
abbreviation $d:=d(C,A)$, we have

$$\cos(\alpha) \ge \frac{{\cal C}^2
-\frac{\varepsilon^2}{2}}{{\cal D}^2}= 1-
d^2\frac{5-d^4}{2(1+d^2)}.$$ So $$\sin^2(\alpha) \le
d^2\left[ \frac{5-d^4}{1+d^2} - d^2\left( \frac
{5-d^4}{2(1+d^2)} \right)^2 \right].$$ Notice that due
to the restriction $d<0.15$ we have $\cos(\alpha)>
\frac{\sqrt{3}}{2}$, so that $|\alpha|< \frac{\pi}{6}
< 0.8$. So
$$|\alpha| \le \frac{8}{7} |\sin(\alpha)| \le
\frac{8}{7}d \cdot \sqrt{ \frac{5-d^4}{1+d^2} -
d^2\left( \frac {5-d^4}{2(1+d^2)} \right)^2} \le 3d
.$$\Kasten\\[0.2cm]
The above corollary estimates the angle of two vectors based at
$\tilde{C}=\expi_A(C)$. Now we will estimate the angle of certain
vectors based at $0\in T_AM$.
\begin{Cor}\label{CorA3}$$ \measuredangle(\tbc -0, _{CA}\parap
\dot{CB}) < 4 d(C,A).$$  \end{Cor} \pf Since $\tb -\tc=\tbc -0$ we
just have to estimate
$$\measuredangle\big({_{CA}\parap}
\dot{CB}\;,\;(\expi_A)_* \dot{CB}\big
)\;\;=\;\;\measuredangle \big({_{f_0}\parap}
\dot{\sigma}(0)\;,\; \dtsigma(0)\big)$$ and apply the
triangle inequality together with Corollary \ref{CorA2}.\\
Denoting by $L$ the distance from $_{f_0}\parap \dot{\sigma}(0)
\in T_AM$ to the line spanned by $ \dot{\tilde{\sigma}}(0)$, we
have $L \le | _{f_0} \parap \dsigma(0) -  \dtsigma(0)| \le
\frac{1}{5}d(C,A)d(C,B)$, where we used Lemma \ref{LemA3} in the
second inequality. Hence
$$\frac{7}{8} \cdot \measuredangle
\big({_{f_0}\parap} \dot{\sigma}(0)\;,\;
\dtsigma(0)\big) <\sin
\left(\measuredangle({_{f_0}\parap}
\,\dot{\sigma}(0)\;,\; \dtsigma(0)) \right) =
\frac{L}{|\dsigma(0)|} \le \frac{1}{5}d(C,A)$$  where
we used $\frac{\sin(x)}{x}>\frac{7}{8}$ for $x\in
[0,0.8]$.\\
Combining this with Corollary \ref{CorA2} gives
$$\measuredangle \left(_{CA}\parap \dot{CB},\tb- \tc
\right)< 3d(A,C) + \frac{8}{35}d(A,C) < 4d(A,C).$$
\Kasten\\[0.2cm]
We conclude this appendix by deriving the estimate need in Claim 3
of Proposition  \ref{Prop31}.
\begin{Cor}\label{CorA4}  $$|\tb -\tc| < \frac{11}{10} d(C,B).$$
\end{Cor}
\pf This follows easily from Lemma \ref{LemA2} since
$$|\tb-\tc| \le \int_0^1 |\dtsigma(s)|ds   \le
\frac{11}{10}d(C,B).$$ \Kasten
\section{An upper bound for $\alpha$ using the curve
$c$}

Here we will prove Proposition  \ref{Prop41}, namely the estimate
 $$\left|\exp_{N_g}^{-1} C- _{\pii}^{\perp} \parap
(\exp_{N_g}^{-1} A) \right| \le L(\gamma) \frac{3150 \e}{f(r)}.$$
To do so we will use the fact that $N$ is $C^1$-close to $N_g$,
see Lemma \ref{LemB3}.\\
In addition to the notation introduced in Section \ref{s5} to
state the proposition, we will use the following.\\

\noindent We will denote by $\pi_c(t)$ the curve $\png \circ
c(t)$, so $\pi_c$ is just a reparametrization of
$\pi$.\\
We will use $\exp$ as a short-hand notation for the
normal exponential map $\exp_{N_g}:(\nu N_g)_{1}
\rightarrow \tnng$. Therefore
$\tilde{c}(t):=\expi(c(t))$ will be a section of $\nu
N_g$ along $\pi_c$.\\
The image under $\exp_*$ of the Ehresmann connection
corresponding to $\np$ will be the subbundle
$\text{LC}^g$ of $TM|_{\exp_{N_g}(\nu N_g)_1}$.\\
To simplify notation we will denote by $\prt \horg$
the projection of $\dgamma(t)\in T_{\gamma(t)}M$ onto
$\vertg_{\gamma(t)}$ \textit{along}
$\horg_{\gamma(t)}$ . We will also use $\prt \ahg$ and
$\prt LC^g$ to denote projections  onto
$\avg_{\gamma(t)}$  \textit{along}  $\ahg_{\gamma(t)}$
and $LC_{\gamma(t)}^g$ respectively.\\[0.4cm]
Our strategy will be to bound above $|\npdt
\tilde{c}(t)|=|\expi_*({pr}_{\dot{c}(t)} \lc)|$ (see Lemma
\ref{LemB3}) using $$TN \approx \horg \approx \ahg \approx \lc.$$
Integration along $\pi_c$ will deliver
the desired estimate.\\[0.4cm]

\noindent The estimates to make precise $TN \approx \horg$ and $ \horg
\approx \ahg$ were derived in [We]. In the next two lemmata we
will do the same for $\ahg \approx \lc$.
\begin{Lem}\label{LemB1} If $L<0.08$ and $p$ is a point in
$\partial \tnngl$, then $$d(\ahg_p, LC^g_p) \le \arcsin
\left(\frac{9}{5}L \right).$$ \end{Lem} \pf It is enough to show
that, if $Y\in \lc_p$ is a
unit vector, $|\pr_Y \ahg| \le \frac{9}{5}L$.\\
Let $\beta(s)$ be a curve tangent to the distribution $\lc$ \st
$\beta(0)=p$,$\,\dot{\beta}(0)=Y$. Then $\expi(\beta(s))=L\xi(s)$
for a unit length parallel section $\xi$ of $\nu N_g$ along the
curve $\gamma(s):=\png(\beta(s))$. If we denote by $K$ the $N_g$
\jf arising from the variation $f_s(t)=\exp(t\xi(s))$, then
clearly $K(L)=Y$  and $K(0)=\dgamma(0)$.\\We claim that $\xi$ is a
strong \jf (see the remark in Section \ref{s32}): we have
$\frac{\partial}{\partial t} |_0 f_s(t)= \xi(s)$, so
$$ K'(0)= \ndt \Bigl|_0 \frac{\partial}{\partial s}
\Big|_0 f_s(t)= \nds \Big|_0 \xi(s)= \npds\Big|_0
\xi(s)-A_{\xi(0)}\dgamma(0)= -A_{\xi(0)}K(0).$$ The claim follows
since $\xi(0)=\dgamma_p(0)$, where $\gamma_p$ denotes the unique
geodesic parametrized to arc-length connecting $\png(p)$ to
$p$.\\Now let us denote by $J$ the $N_g$ \jf along $\gamma_p$
vanishing at $0$ \st $J(L)=\pr_Y \ahg \in \avg_p$. By Lemma
\ref{Lem29}, using the fact that $Y$ is a unit vector,  we have
$$|\pr_Y \ahg|^2=\langle \pr_Y \ahg, Y \rangle =
|\langle J(L),K(L) \rangle| \le \frac{9}{5}L\cdot
|\pr_Y \ahg|$$ and we are done.\Kasten\\[0.2cm]
\begin{Lem}\label{LemB2} Let $L<0.08$. For any point $p$ in $
\partial \exp_{N_g}(\nu N_g)_L$ the projections $T_pM
\rightarrow \avg_p$ along $\ahg_p$ and $LC_p^g$ differ
at most by $2L$ in the operator norm. \end{Lem}
\pf Let $\phi: \ahg_p \rightarrow \avg_p$ be the
linear map whose graph is $\lc_p$.
Let $X\in T_pM$ a unit vector and write $X=X_{ah}+X_{av}$ for the decomposition of $X$ into almost horizontal
and almost vertical vectors. Then $X=(X_{ah} +\phi(X_{ah}))
+(X_{av}-\phi(X_{ah}))$ is the decomposition with respect to the subspaces  $\lc_p$ and $\avg_p$. The difference of the
two projections onto $\avg_p$ maps $X$ to $\phi(X_{ah})$. Now
$$ |\phi(X_{ah})| \le |\phi|_{op} \le \tan (d(\ahg_p, \lc_p))\le \frac{\frac{9}{5}L}{\sqrt{1 -
(\frac{9}{5}L)^2}}<2L,$$
where we used [We, Cor. A.5] in the second inequality and Lemma \ref{LemB1} in the third one. \Kasten\\[0.2cm]

\noindent Now we are ready to bound the covariant derivative of
$\tilde{c}(t)$:
\begin{Lem}\label{LemB3} For all $t$   $$ \left| \npdt \tilde{c}(t)
\right| \le 2702 \e.$$ \end{Lem} \pf Let $\widehat{\npdt
\tilde{c}(t)}$ denote $\npdt \tilde{c}(t) \in  \nu_{\pi_c(t)}N_g$
but considered as an element of $T_{\tilde{c}(t)}(
\nu_{\pi_c(t)}N_g)$. First notice that, by definition,
$\widehat{\npdt \tilde{c}(t)}$ is the image of
$\dot{\tilde{c}}(t)$ under the projection $T_{{\tilde{c}}(t)}(\nu
N_g) \rightarrow T_{{\tilde{c}}(t)}(\nu_{\pi_c(t)} N_g)$ along the
Ehresmann connection on $\nu N_g$ corresponding to $\np$.
Therefore, since $\exp_*$ maps the this Ehresmann connection to
$LC^g$ and tangent spaces to the fibers of $\nu N_g$ to $\avg$, we
have
$$\exp_*\left(\widehat{\npdt \tilde{c}(t)}\right)=
\pr_{\dot{c}(t)} \lc.$$ Notice that here $\exp_*$
denotes $d_{\tilde{c}(t)}\exp_{N_g}$.\\[0.2cm]
The fact that $N$ is $C^1$-close the $N_g$ (see Theorem \ref{Thmk}),
since $\dot{c}(t) \in T_{c(t)}N$ implies that  $\measuredangle
(\dot{c}(t) ,
 \horg_{c(t)}) \le 2500 \e$. By [We, Prop. 3.7]
$d(\horg_{c(t)}, \ahg_{c(t)})
\le \frac{\e}{4}$ since $d(c(t),N_g) \le 100 \e$.\\
Therefore $\measuredangle (\dot{c}(t) , \ahg_{c(t)}) \le 2501 \e $
and $|\pr_{\dot{c}(t)} \ahg|\le \sin( 2501 \e) \le 2501 \e.$ \\On
the other hand, by Lemma \ref{LemB2}, $|\pr_{\dot{c}(t)} \ahg
-\pr_{\dot{c}(t)} \lc| \le 200 \e$. The triangle inequality
therefore gives $|\pr_{\dot{c}(t)} \lc|\le \sin( 2701 \e)$.
Therefore, using Lemma \ref{LemA2} and $\e < \frac{1}{20000}$,
$$|\expi_* (\pr_{\dot{c}(t)} \lc)| \le \frac{1}{1
-\frac{\e}{5}} \left|\pr_{\dot{c}(t)} \lc\right|
\le 2702 \e.$$\Kasten\\[0.2cm]
Lemma \ref{LemB3} allows us to bound $\left| \expi C -
_{\pi^b}^{\perp}
\parap \left(\expi A \right) \right|$ in terms of $L(c)$. However
we want a bound in terms of $L(\gamma)$, so
now we will now compare the lengths the two curves.\\
Recall that $f(x) =\cos(x)- \frac{3}{2}\sin(x)$ and $r:=100\e +
\frac{L(\gamma)}{2}$. Notice also that $r<0.08$ due to our
restrictions on $\e$ and $d(C,A)$.
\begin{Lem}\label{LemB4}$$ L(c)\le \frac{1+ 3200
\e}{f(r)}L(\gamma).$$ \end{Lem} \pf Since
$\vg_*({\png}_*\dot{c}(t)) = \dot{c}(t)$, by Proposition
\ref{Prop21} we have $|\dot{c}(t) - \para ({\png}_*\dot{c}(t))|\le
3200 \e|{\png}_*\dot{c}(t)|$, so
$$|\dot{c}(t)|\le (1 + 3200 \e) |
{\png}_*\dot{c}(t)|.$$ Since $L(\png \circ c)=L(\pi)$, from this
follows $L(c) \le (1+3200 \e) L(\pi)$. By [We, Lemma 3.3] we have
$f(r)L(\pi) \le L(\gamma)$ and
we are done. \Kasten\\[0.2cm]
\textit{ Proof of Proposition  \ref{Prop41}: } We have
\begin{equation*} \begin{split}
\left| \expi C - _{\pi^b}^{\perp} \parap \left(\expi A
\right) \right|=&
\left| \int_0^{L(c)} \frac{d}{dt}\;
_{\pi_c^b}^{\perp}  \parap \tilde{c}(t) \text{dt}\right|
\\
=& \left|{ \int_0^{L(c)}} {_{\pi_c^b}^{\perp} \parap}
\npdt  \tilde{c}(t) \text{dt}\right| \\
\le&  2702 \e L(c) \\
\le&   2702 \e \frac{(1+3200 \e)}{f(r)} L(\gamma)
\end{split} \end{equation*} where we used Lemma \ref{LemB3}
and lemma \ref{LemB4} in the last two inequalities. The
proposition follows using the bound $\e < \frac{1}{20000}$.\Kasten

\section{A lower bound for $\alpha$ using the curve
$\gamma$}

Here we will prove Proposition  \ref{Prop42}, i.e. the estimate
$$\left|\expi_{N_g}(C) -  _{\pi^b}^{\perp} \parap
\expi_{N_g}(A)\right| \ge L(\gamma)\left[ \frac{99}{100}\sin
\left(\alpha -\frac{\epsilon}{4} \right) - 500 \e - 3r -
\frac{8}{3}L(\gamma)\left(r+\frac{r + \frac{3}{2}}{f(r)}
\right)\right].$$ We will use the fact that $N_g$ has bounded
second fundamental form (see the first statement of Lemma
\ref{LemC3}) and that $\gamma$ is a geodesic (see the second
statement of the same
Lemma).\\
 We will use the notation introduced in Section \ref{s5}
and at the beginning of Appendix B. Recall that
$\tilde{\gamma}(t):=\expi_{N_g}(\gamma(t))$ is a
section of $\nu N_g$ along $\pi$.\\[0.4cm]
First we will set a lower bound on the initial derivative of
$\tilde{\gamma}$.
\begin{Lem}\label{LemC1} We have $$ \left| \npdto \tgamma(0)
\right| \ge
\frac{99}{100}\left[\sin\left(\alpha-\frac{\epsilon}{4}\right)
-200 \e\right].$$ \end{Lem} \pf Analogously to the proof of Lemma
\ref{LemB3} we have $\exp_*(\widehat{\npdt \tgamma(0)})= \pro
\lc$, where $\widehat{\npdt \tgamma(0)}$ is an element of
$T_{\tilde{\gamma}(0)} \nu_{\pi(0)}N_g$.\\
By [We, Prop 3.7] we have $d(\horg_C\;,\;\ahg_C) \le
\frac{\e}{4}$. So
$$\measuredangle \left(\dgamma (0) \;,\; \ahg_{C}
\right) \ge \measuredangle \left(\dgamma(0)\; ,\;
\horg_{C} \right) - d(\horg_C\;,\; \ahg_C) \ge \alpha
- \frac{\e}{4}.$$
Therefore $|\pro \ahg| \ge \sin(\alpha -
\frac{\e}{4}).$\\
On the other hand, by Lemma \ref{LemB2}, $|\pro \ahg - \pro \lc|
\le 200 \e$. The inverse triangle inequality gives $$| \pro \lc|
\ge \sin(\alpha - \frac{\e}{4}) - 200 \e.$$ Applying $\expi_*$, by
Lemma \ref{LemA2} we have $|\expi_*( \pro \lc)| \ge
\frac{1}{1+\frac{\e}{5}} | \pro \lc|$, and since $
\frac{1}{1+\frac{\e}{5}}  \ge
\frac{99}{100}$ we are done. \Kasten\\[0.2cm]
Our next goal is to show that $\tgamma(t)$ ``grows at a nearly
constant rate''. This will be achieved in Corollary \ref{CorC3}.
Together with Lemma \ref{LemC1} and integration along $\pi$ this
will deliver the estimate of
Proposition  \ref{Prop42}.\\[0.4cm]
The next two lemmas will be used to prove Corollary \ref{CorC1},
where we will show that $\npdto \tgamma(0)$ and $\expi_* \circ
_{\gammai}\parap \circ \exp_*(\widehat{\npdto \tgamma(t)})$ - i.e.
the parallel translate of $\npdto \tgamma(t)$ ``along $\gamma$'' -
are close for all $t$. Here $\widehat{\npdto \tgamma(t)}$ denotes
the vector ${\npdto \tgamma(t)}$ regarded as an element of
$T_{\tgamma(t)}(\nu _{\pi(t)}N_g)$. To this aim we show that
$$\text{pr}_{\dgamma(0)} \lc \approx
\text{pr}_{\dgamma(0)} \horg \approx \text{pr}_{\dgamma(t)} \horg
\approx \text{pr}_{\dgamma(t)} \lc,$$ where we identify tangent
spaces by parallel translation along $\gamma$. The crucial step is
the second ``$\approx$'', where we use that fact that $\gamma$ is
a geodesic. Applying $\expi_*$ will easily imply Corollary
\ref{CorC1} since $\exp^{-1}_*( \text{pr}_{\dgamma(t)}
\lc)=\widehat{\npdt \tilde{\gamma}(t)}$.
\begin{Lem}\label{LemC2} For any $L<1$ and any point $p\in
\exp_{N_g}(\nu N_g)_L$ the  orthogonal projections $T_pM
\rightarrow \ahg_p$  and $T_pM \rightarrow \horg_p$  differ at
most by $\frac{L^2}{5}$ in the operator norm. \end{Lem} \pf This
follows immediately from [We, Prop.
3.7].\Kasten\\[0.2cm]
\begin{Lem}\label{LemC3}  For all $t$ $$d(\vertg_C,_{\gamma^b}
\parap \vertg_{\gamma(t)}) \le \arcsin \left[t \left(r+
\frac{r+\frac{3}{2}}{f(r)} \right) \right].$$
Furthermore,
 $$ \big|\pro \horg - _{\gammai} \parap \prt \horg
\big| \le t \left(r+ \frac{r+\frac{3}{2}}{f(r)}
\right).$$\end{Lem}
\pf We first want to estimate $d(\vertg_C\;,\;
_{\gamma^b} \parap \vertg_{\gamma(t)})$.
Let $v\in \nu_C N_g$ be a normal unit vector.\\
First of all, for the $\nabla$ and $\np$ parallel
translations along $\pi$ from $C$ to $\pi(t)$ we have
$$\left| _{\pi}\parap v -  _{\pi}^{\perp} \parap v
\right| \le \frac{3}{2}L(\pi|_{[0,t]}) \le
\frac{3}{2}\frac{t}{f(r)}.$$ The first inequality
follows from a simple computation involving the second
fundamental form of $N_g$, which is bounded in norm by
$\frac{3}{2}$ (see [We, Cor. 3.2]). The second
inequality is due to $f(r) L(\pi|_{[0,t]})\le
L(\gamma|_{[0,t]})$, which follows from [We, Lemma
3.3].\\
Secondly, denoting by $\tau_t$ the unit speed geodesic
from $\pi(t)$ to $\gamma(t)$, we have
$$\Big| {_{\tau_0}\parap v}\; -\;  _{\gamma^b}
\parap \circ\;  _{\tau_t} \;\para \circ\;  _{\pi}
\parap v \Big| \le rt \left(1+ \frac{1}{f(r)}
\right).$$ Indeed, the above expression just measures the holonomy
as one goes once around the polygonal loop given by the geodesics
$\tau_0^b, \pi,\tau_t$ and $\gamma^b$. Using the bounds on
curvature we know that this is bounded by the area of a surface
spanned by the polygon (see [BK, 6.2.1]). The estimate given above
surely holds since $L(\tau_t),L(\tau_0) \le r$, $
L(\gamma|_{[0,t]})=t$ and, as we just saw, $
L(\pi|_{[0,t]}) \le \frac{t}{f(r)}$.\\
Together this gives  \begin{equation*} \begin{split}
\Big| {_{\tau_0}\parap v}\; -\;  _{\gamma^b} \parap
\circ \; _{\tau_t} \parap \circ\;  _{\pi}^{\perp}
\parap v \Big|  \le& \Big| {_{\tau_0}\parap v}\; -\;
_{\gamma^b} \parap \circ \; _{\tau_t}\parap \circ
\; _{\pi} \parap v \Big| \;+\; \Big|
{_{\gamma^b}\parap} \;  \circ\;  _{\tau_t} \parap
\circ\; \left[_{\pi}\parap v -  _{\pi}^{\perp} \parap
v \right]\Big|\\
\le& t\left( r + \frac{r +\frac{3}{2}}{f(r)} \right).
\end{split} \end{equation*}
So we obtain a bound on the distance from $ _{\tau_0}\parap v \in
\vertg_C$ to a unit vector in $ _{\gamma^b}\parap
\vertg_{\gamma(t)}$. This delivers the first statement of the
lemma.\ The second statement follows using [We, Prop. A.4], since
$ _{\gamma^b}\parap \pr_{\dgamma(t)} \horg = {\pr_{\dgamma(0)}}
(_{\gamma^b} \parap \horg_{\gamma(t)})$ because $\gamma$
is a geodesic.\Kasten\\[0.2cm]
\begin{Cor}\label{CorC1}  For all $t$ $$\left|\expi_* \circ
_{\gammai}\parap \circ \exp_*  \left(\widehat{\npdto \tgamma(t)}
\right) -  \npdto \tgamma(0)\right|
 \le \frac{51}{50}\left[2.1(100\e+r) +t\left(r+
\frac{r+\frac{3}{2}}{f(r)}\right)\right].$$ \end{Cor} \pf From
Lemma \ref{LemC2} and Lemma  \ref{LemB2} we have for all $t$

\begin{equation*}\begin{split} \big|\prt \horg -\prt
LC^g \big|  & \le   \big|\prt \horg -\prt \ahg \big|+
\big|\prt \ahg -\prt LC^g \big| \\ & \le
\frac{r^2}{5} + 2r \\  & \le  2.1 r .\end{split}
\end{equation*}
For $t=0$, since $d(C, N_g) < 100\e$, we have the better estimate
$$\big|\pro \horg -\pro LC^g\big|  \le 210 \e.$$ Combining this
with the second statement of Lemma \ref{LemC3} gives $$\left| \pro
LC^g - _{\gammai} \parap \prt LC^g\right| \le 2.1(100 \e +r) +
t\left(r+
\frac{r+\frac{3}{2}}{f(r)}\right).$$\\
Recall that $\prt LC^g = \exp_* (\widehat{\npdt \tgamma(t)})$, as
in the proof of Lemma \ref{LemB3}. Also, for any vector $X\in
T_CM$ we have $|\expi_* X| \le \frac{|X|}{1-\frac{\epsilon}{5}}$
by Lemma \ref{LemA2}. So applying $(\expi)_*$ to $ \pro LC^g -
_{\gammai}
\parap \prt LC^g$ we get
 $$  \left| \npdto \tgamma(0)-
\expi_* \circ _{\gammai}\parap \circ \exp_*
\left(\widehat{\npdto \tgamma(t)} \right) \right| \le
\left[2.1(100\e+r) +t\left(r+
\frac{r+\frac{3}{2}}{f(r)}\right)\right]
\frac{1}{1-\frac{\epsilon}{5}}.$$ \Kasten\\[0.2cm]
Now let $\xi$ be a unit vector in $\nu_{\pi(t)}N_g$. Denote by
$\hat{\xi}$ the same vector thought of as an element of
$T_{\tgamma(t)}(\nu_{\pi(t)}N_g)$.\\ In the next two lemmas we
want to show that $ _{\pii}^{\perp} \parap {\xi}$ and $\expi_*
\circ _{\gammai}\parap \circ \exp_* \hat{\xi} \in T_CM$  are close
to each other, i.e. that under the identification by $\exp$ the
$\nabla^{\perp}$-parallel translation along $\pi$ and the
$\nabla$-parallel translation along $\gamma$ do not differ too
much. Here we also make use of the fact that $N$ has bounded
second fundamental form (see Lemma \ref{LemC5}). In Corollary
\ref{CorC2} we will apply this to the vector $\npdto \tgamma(t)$.
\begin{Lem}\label{LemC4} Denoting by $\tau_t$ the unit speed
geodesic from $\pi(t)$ to $\gamma(t)$, $$\left|
_{\tau_0^b}\parap\circ  _{\gammai}\parap \circ _{\tau_t}\parap
{\xi} \;-\; \expi_* \circ _{\gammai}\parap \circ \exp_* \hat{\xi}
\right| < \frac{r^2}{2}.$$ \end{Lem} \pf First let us notice that
applying Lemma \ref{LemA2} three times we get
\begin{equation*} \begin{split}
\left|  _{\tau_0^b} \parap  \left[  _{\gamma^b}\,
\para \exp_* \hat{\xi} \right] - \expi_* \left[
 _{\gamma^b} \parap \exp_* \hat{\xi} \right] \right|
\le& \frac{r^2}{5} \left|  \expi_* \left[
 _{\gamma^b} \parap \exp_* \hat{\xi} \right]
\right|\\
\le &  \frac{r^2}{5} \frac{1}{1-\frac{r^2}{5}} \left|
_{\gamma^b}  \parap \exp_* \hat{\xi} \right|\\
\le &  \frac{r^2}{5} \frac{ 1 +\frac{r^2}{5}}{1-\frac{r^2}{5}}.
\end{split}
\end{equation*}
Therefore, applying Lemma \ref{LemA2} to $\xi$,  the left hand
side of the statement of this lemma is bounded above by
\begin{equation*} \begin{split}  \left|  _{\tau_0^b}
\parap \circ\;  _{\gamma^b} \parap \Big[  _{\tau_t}
\parap \xi \Big] -
 _{\tau_0^b} \parap \circ\;  _{\gamma^b} \parap
\Big[ \exp_* \hat{ \xi} \right] \Big|
+& \left|  _{\tau_0^b}\; \para \circ \left[
_{\gamma^b}\; \para \exp_* \hat{ \xi} \right]
- \expi_*  \left[  _{\gamma^b} \parap \exp_* \hat{
\xi} \right] \right| \\
\le & \frac{r^2}{5} +  \frac{r^2}{5} \frac{ 1
+\frac{r^2}{5}}{1-\frac{r^2}{5}} \\
\le& r^2\frac{2}{5(1- \frac{r^2}{5})}. \end{split}
\end{equation*} \Kasten\\[0.2cm]
\begin{Lem}\label{LemC5} $$\left|\expi_* \circ _{\gammai}\parap
\circ \exp_* \hat{\xi}\; -\;  _{\pii}^{\perp} \parap \;{\xi}
\right| \le \frac{r^2}{2} + t \left(r+ \frac{r +
\frac{3}{2}}{f(r)} \right).$$ \end{Lem} \pf The left hand side in
the statement of the lemma is bounded above by
\begin{equation*} \begin{split}
& \left|  \expi_* \circ\;  _{\gamma^b} \parap
\,\;\circ\; \exp_* \hat{\xi} \;- \;
 _{\tau_0^b}\parap \circ\;  _{\gamma^b} \parap
\circ\;  _{\tau_t} \parap \xi \right|\\
+ & \left| _{\tau_0^b} \parap \circ\;  _{\gamma^b}
\parap \circ\;  _{\tau_t} \parap \xi \;-\;
 _{\pi^b} \parap \xi \right| \\
+& \left|  _{\pi^b}\parap \xi \;-\;
_{\pi^b}^{\perp} \parap\, \xi \right|\\
\le& \frac{r^2}{2} + rt\left( 1+ \frac{1}{f(r)}
\right) + \frac{3}{2} \frac{t}{f(r)} .
 \end{split} \end{equation*} The first term is
estimated by Lemma \ref{LemC4}.\\
The second one is just the holonomy as once goes around the loop
given by $\tau_t, \gamma^b, \tau_0^b$ and $\pi$, which was bounded
above in the proof of
Lemma \ref{LemC3}.\\
The third and last term is estimated in the proof of
Lemma \ref{LemC3} as well. \Kasten\\[0.2cm]
\begin{Cor}\label{CorC2} The section $\tgamma$ satisfies
$$\left|\expi_* \circ _{\gammai}\parap \circ  \exp_*
\left(\widehat{\npdto \tgamma(t)}\right) -\, _{\pii}^{\perp}
\parap  {\npdto \tgamma(t)} \right| \le \frac{2}{3}r^2   +
\frac{4}{3}t\left( r + \frac{r + \frac{3}{2}}{f(r)} \right).$$
\end{Cor} \pf We apply Lemma \ref{LemC5} to $ \widehat{\npdto
\tgamma(t)}$, where now we have to take into consideration the
length of $ \widehat{\npdto
\tgamma(t)}$ in our estimate.\\
We have $| \widehat{\npdto \tgamma(t)}| =|\expi_*(\prt LC^g)|\le
\frac{1}{1- \frac{r^2}{5}} |\prt LC^g|$ by Lemma \ref{LemA2}, and
$$|\prt LC^g| \le |\prt LC^g - \prt \ahg| + |\prt \ahg| \le 2r + 1
$$ by Lemma \ref{LemB2}. Since $\frac{2r+1}{1-\frac{r^2}{5}} \le
\frac{4}{3}$ for
$r\le 0.08$ we are done. \Kasten\\[0.2cm]
Now Corollary \ref{CorC1} and  Corollary \ref{CorC2} immediately
imply that $\tilde{\gamma}(t)$ ``grows at a nearly constant
rate'':
\begin{Cor}\label{CorC3} The section $\tgamma$ satisfies
$$\left|\npdt \tgamma(0) \;-  \;_{\pii}^{\perp}\parap
 {\npdto \tgamma(t)} \right | \le 3(100 \e +r) +
\frac{8}{3}t \left(r+\frac{r + \frac{3}{2}}{f(r)} \right) .$$
\end{Cor} \textit{Proof of Proposition  \ref{Prop42}: } The
estimate of Proposition  \ref{Prop42} follows from
\begin{equation*}
\begin{split}  \left|\expi(C) -  _{\pi^b}^{\perp}\parap \expi(A) \right|=& \left|{\int_0^{L(\gamma)}} {_{\pi^b}^{\perp} \parap}  \npdt \tgamma(t)dt \right|\\  \ge & \left|{\int_0^{L(\gamma)}}  \npdt \tgamma(0)dt \right|-\left|{\int_0^{L(\gamma)}} \left(\npdt \tgamma(0)-  _{\pi^b}^{\perp} \parap \npdt \tgamma(t) \right) dt\right|\\ \ge& L(\gamma)\cdot \left| \npdt \tgamma(0) \right| - {\int_0^{L(\gamma)}} \left|\npdt \tgamma(0)- _{\pi^b}^{\perp} \parap \npdt \tgamma(t) \right| dt \end{split} \end{equation*}
using Lemma \ref{LemC1} and Corollary \ref{CorC3}. \Kasten

\end{document}